\RequirePackage{ifthen}
\newboolean{DCG}
\setboolean{DCG}{false}

\ifthenelse {\boolean{DCG}}
{
\documentclass[smallextended]{svjour3}
\smartqed
} {
\documentclass[11pt]{article}
\usepackage[margin=1.3in]{geometry}
}

\usepackage{amsmath,amssymb}
\usepackage{graphicx}
\usepackage[margin=1.3in]{geometry}
\usepackage{color}
\usepackage{amsfonts}
\usepackage{hyperref}
\usepackage{soul}
\usepackage{cancel}
\usepackage{enumitem}

\usepackage{algorithm}
\usepackage[noend]{algpseudocode}

\floatname{algorithm}{}

\usepackage{xspace}
\newcommand{\algstyle}[1]{\textsf{#1}\xspace}
\newcommand{\oracle}{\algstyle{oracle}}
\newcommand{\Oracle}{\algstyle{Oracle}}
\newcommand{\basic}{\algstyle{basic algorithm}}
\newcommand{\Basic}{\algstyle{Basic algorithm}}
\newcommand{\preprocessing}{\algstyle{preprocessing algorithm}}
\newcommand{\Preprocessing}{\algstyle{Preprocessing algorithm}}
\newcommand{\scaling}{\algstyle{scaling algorithm}}
\newcommand{\Scaling}{\algstyle{Scaling algorithm}}
\newcommand{\preprocessingscaling}{\algstyle{preprocessing \& scaling algorithm}}
\newcommand{\Preprocessingscaling}{\algstyle{Preprocessing \& scaling algorithm}}
\newcommand{\iterative}{\algstyle{iterative algorithm}}
\newcommand{\Iterative}{\algstyle{Iterative algorithm}}

\newcommand{\ceil}[1]{\left\lceil#1\right\rceil}
\newcommand{\floor}[1]{\left\lfloor#1\right\rfloor}

\newcommand{\abs}[1]{\left\lvert#1\right\rvert}
\newcommand{\norm}[1]{\left\lVert#1\right\rVert_2}
\newcommand{\norminf}[1]{\left\lVert#1\right\rVert_\infty}
\newcommand{\normone}[1]{\left\lVert#1\right\rVert_1}

\newcommand{\B}{\mathcal B}
\newcommand{\E}{\mathcal E}
\renewcommand{\H}{\mathcal H}

\newcommand{\R}{\mathbb R}
\newcommand{\Z}{\mathbb Z}
\newcommand{\Q}{\mathbb{Q}}

\DeclareMathOperator{\argmax}{argmax}

\DeclareMathOperator{\sign}{sign}
\DeclareMathOperator{\size}{size}

\ifthenelse {\boolean{DCG}}
{
\newenvironment{prf}[1][]
{\begin{proof}}
{\qed \end{proof}}

\newenvironment{prfc}[1][]
{\begin{proof}[#1]}
{\qed \end{proof}}

\journalname{Discrete \& Computational Geometry}

} {
\usepackage{amsthm}

\newtheorem{theorem}{Theorem}
\newtheorem{lemma}{Lemma}
\newtheorem{corollary}{Corollary}
\newtheorem{remark}{Remark}
\newtheorem{proposition}{Proposition}

\newenvironment{prf}[1][]
{\begin{proof}}
{\end{proof}}

\newenvironment{prfc}[1][]
{\begin{proof}[Proof #1]}
{\end{proof}}

}

\newtheorem{observation}{Observation}

\begin{document}

\title{Short simplex paths in lattice polytopes}
%\title{A strongly polynomial simplex method for lattice polytopes in small boxes}
%\title{The simplex method is polynomial for lattice polytopes}
%\title{The simplex method is strongly polynomial for simple lattice polytopes}
%\title{An optimal simplex algorithm for lattice polytopes}

\ifthenelse {\boolean{DCG}}
{
% For DCG begin
\titlerunning{Short simplex paths in lattice polytopes}
\authorrunning{Alberto Del Pia, Carla Michini}

\author{Alberto~Del~Pia \and Carla~Michini}
\institute{Alberto~Del~Pia \at
              Department of Industrial and Systems Engineering 
              \& Wisconsin Institute for Discovery \\
              University of Wisconsin-Madison, Madison, WI, USA \\
              \email{delpia@wisc.edu}
              \and
              Carla~Michini \at
              Department of Industrial and Systems Engineering \\
              University of Wisconsin-Madison, Madison, WI, USA \\
              \email{michini@wisc.edu}}
% For DCG end
}
{
% For OO begin
\author{Alberto~Del~Pia
\thanks{Department of Industrial and Systems Engineering \& Wisconsin Institute for Discovery,
             University of Wisconsin-Madison, Madison, WI, USA.
             E-mail: {\tt delpia@wisc.edu}.}
\and
Carla~Michini
\thanks{Department of Industrial and Systems Engineering,
             University of Wisconsin-Madison, Madison, WI, USA.
             E-mail: {\tt michini@wisc.edu}.}
             }
% For OO end
}

\date{\today}
%\date{March 09, 2020}

\maketitle

\begin{abstract}
The goal of this paper is to design a simplex algorithm for linear programs on lattice polytopes that traces `short' simplex paths from any given vertex to an optimal one.
We consider a lattice polytope $P$ contained in $[0,k]^n$ and defined via $m$ linear inequalities. 
Our first contribution is a simplex algorithm that reaches an optimal vertex by tracing a path along the edges of $P$ of length in $O(n^4 k\log(nk))$.
The length of this path is independent from $m$ and it is the best possible up to a polynomial function.
In fact, it is only polynomially far from the worst-case diameter, which roughly grows as a linear function in $n$ and $k$.
%can grow essentially linearly in both $n$ and $k$.

Motivated by the fact that most known lattice polytopes are defined via $0,\pm 1$ constraint matrices, our second contribution is an iterative algorithm which exploits the largest absolute value $\alpha$ of the entries in the constraint matrix.
We show that the length of the simplex path generated by the iterative algorithm is in $O(n^2k \log(nk\alpha))$.
In particular, if $\alpha$ is bounded by a polynomial in $n, k$, then the length of the simplex path 
is in $O(n^2k \log(nk))$.

For both algorithms, the number of arithmetic operations needed to compute the next vertex in the path is polynomial in $n$, $m$ and $\log k$.
If $k$ is polynomially bounded by $n$ and $m$, the algorithm runs in strongly polynomial time.

\ifthenelse {\boolean{DCG}}
{
% For DCG begin
% For DCG end
} {}
\end{abstract}

\ifthenelse {\boolean{DCG}}
{}{
% For OO begin
}

% ===============================================
% ===============================================

\section{Introduction}
\label{sec intro}

Linear programming (LP) is one of the most fundamental types of optimization models.
%a pillar for computation and theory of mathematical optimization and a powerful modeling tool in many applications \note{strange sentence}.
In an LP problem, we are given a polyhedron $P \subseteq \R^n$ and a cost vector $c \in \Z^n$, and we wish to solve the optimization problem
%
%
%
%We consider the LP problem
%In this paper we are interested in giving a simplex algorithm for the LP problem
\begin{align}
\label{prob: main}
\max \{c^\top x \mid x \in P\}.
\end{align}
The polyhedron $P$ is explicitly given via a 
system of 
%, given via an external description consisting of 
linear inequalities, i.e.,
$P = \{x \in \R^n \mid Ax \le b\}$,
where $A \in \Z^{m \times n}$, $b \in \Z^m$.
%, the problem of optimizing a linear function over a polyhedron, 
%Note that, i
If $P$ is nonempty and bounded, problem \eqref{prob: main} admits an optimal solution that is a vertex of $P$.

In this paper, we consider a special class of linear programs \eqref{prob: main} where the feasible region is a \emph{lattice polytope}, i.e.,~a polytope whose vertices have integer coordinates. 
These polytopes are particularly relevant in discrete optimization and integer programming, as they correspond to the convex hull of the feasible solutions to such optimization problems.
In particular, a \emph{$[0,k]$-polytope} in $\R^n$ is defined as a lattice polytope contained in the box $[0,k]^n$.  
%\cnote{can we assume that only the vertices are in the box?}

One of the main algorithms for LP is the {simplex method}.
%The simplex method is one of the main algorithms for LP, and has been selected as one of the most influential algorithms in the $20$th century \cite{DonSul00}.
%
% However, its running time is not always polynomial in the input size \note{NOT SURE}.
%
%it is not known to be polynomial 
%or strongly polynomial 
%\note{refer to a book for size and related complexity stuff.}
%We remark that in this paper we use the standard notion of polynomial-time algorithms in Discrete Optimization, and we refer the reader to Section~2.4 in the book~\cite{SchBookIP} for a thorough introduction.
The simplex method moves from the current vertex to an adjacent one along an edge of the polyhedron $P$, until an optimal vertex is reached or unboundedness is detected, and the selection of the next vertex depends on a {pivoting rule}.
The sequence of vertices generated by the simplex method is called the \emph{simplex path}.
The main objective of this paper is to design a simplex algorithm for $[0,k]$-polytopes that constructs `short' simplex paths from any starting vertex $x^0$.

%The length of the simplex path might be exponential in the input size \note{NOT SURE}, and a fundamental question in linear programming is whether there exists a pivoting rule yielding a polynomial time variant of the simplex method. 
%
%\cnote{A lower bound on the number of vertices of $P$ that the simplex method might visit is the \emph{diameter} of $P$.}
{But how short can a simplex path be?}
A natural lower bound on the length of a simplex path from $x^0$ to $x^*$ is given by the \emph{distance} between these two vertices, which is defined as the minimum length of a path connecting $x^0$ and $x^*$ along the edges of the polyhedron $P$.
The \emph{diameter} of $P$ 
%We recall that the diameter of $P$ 
is the largest {distance} between any two vertices of $P$, 
%where the distance between two vertices is 
and therefore it provides a lower bound on the length of a worst-case simplex path on $P$.
%
%\note{NOT SURE}
%For general linear programs, the known deterministic pivoting rules may visit exponentially many vertices \cite{FriHanZwi11}\cnote{more references}. On the other hand, few special cases admit a polynomial or strongly polynomial pivoting rule, e.g., flow problems \cite{GolHao92,Orl97,GolChe97,Arm98}, the stable set problem on bipartite graphs \cite{IkuNem82,ArmJin96}, Markov decision problems with a fixed discount rate \cite{Ye11} and the linear fractional assignment problem \cite{Kab08}.\cnote{more Mizuno references}
%
%\emph{But how short can a simplex path be for a $[0,k]$-polytope?}
%The diameter of $[0,k]$-polytopes provides a lower bound on the worst-case simplex path length.
It is known that the diameter of $[0,1]$-polytopes in $\R^n$ is at most $n$ \cite{Nad89} and this bound is attained by the hypercube $[0,1]^n$. This upper bound was later generalized to $nk$ for general $[0,k]$-polytopes in $\R^n$ \cite{KleOnn92}, and refined to $n \floor{(k - \frac 1 2)}$ for $k\ge2$ \cite{dPMic16} and to $n k - \ceil{\frac 23 n } - (k-3)$ for $k\ge3$ \cite{DezPou18}.
%
%For general $[0,k]$-polytopes in $\R^n$ Kleinshmidt and Onn \note{names} proved that the diameter is at most $nk$ \cite{KleOnn92}. This upper bound has been later refined to $n \floor{(k - \frac 1 2)}$ for $k\ge2$ \cite{dPMic16} and to $n k - \ceil{\frac 23 n } - (k-3)$ for $k\ge3$ \cite{DezPou18}.
For $k=2$ the bound given in \cite{dPMic16} is tight.
In general, for fixed $k$, the diameter of lattice polytopes can grow linearly with $n$, since there are lattice polytopes, called {primitive zonotopes}, that can have diameter in $\Omega (n)$ \cite{DezManOnn18}.
Viceversa, when $n$ is fixed, the diameter of a $[0,k]$-polytope in $\R^n$ can grow almost linearly with $k$.
In fact, it is known that for $n=2$ there are $[0,k]$-polytopes with diameter in $\Omega (k^{2/3})$ \cite{BalBar91,Thi91,AckZun95}. %Zie95
Moreover, for any fixed $n$, there are primitive zonotopes with diameter %$\Omega()$ \cite{DezPouSuk19}
in $\Omega(k^{\frac {n}{n+1}})$ for $k$ that goes to infinity \cite{DezPouSuk19}.
%Thus, our goal is
%to trace simplex paths on lattice polytopes in $[0,k]^n$ whose length only depends on $k$ and $n$ and is reasonably close to the worst-case diameter.
%to design a simplex algorithm for lattice polytopes in $[0,k]^n$ that constructs simplex paths whose length $(i)$ only depends on $n$ and $k$ and $(ii)$ is reasonably close to the worst-case diameter.
%
%The main objective of this paper is to design an algorithm for $[0,k]$-polytopes that constructs simplex paths whose length is upper bounded by a polynomial of $k$ and $n$.

Can we design a simplex algorithm whose simplex path length is only polynomially far from optimal, meaning that it is upper bounded by a polynomial function of the worst-case diameter?
In this paper, we answer this question in the affermative.
Our first contribution is a \preprocessingscaling
that generates a simplex path of length in $O(n^4 k \log (nk))$ for $[0,k]$-polytopes in $\R^n$.
%\begin{proposition}
%\label{th 1}
%The length of the simplex path generated by the preprocessing \& scaling algorithm is bounded by $nk (4n^3 \log 2 + n(n+2) \log (nk+1) + 2) \in %O(n^4 k \log (nk))$.
%\end{proposition}
%Thus, the length of our simplex path is only polynomially far from optimal, meaning that it is upper bounded by a polynomial function of the worst-case diameter.
The length of the simplex path is indeed polynomially far from optimal, as it is polynomially bounded in $n$ and $k$.
We remark that the upper bound is independent on $m$. This is especially interesting because, even for $[0,1]$-polytopes, $m$ can grow exponentially in $n$ (see, e.g., \cite{SchBookCO}).

Our next objective is that of decreasing the gap between the length $O(n^4 k \log (nk))$ provided by the \preprocessingscaling and the worst-case diameter, for wide classes of $[0,k]$-polytopes.
We focus our attention on $[0,k]$-polytopes with bounded parameter $\alpha$, defined as the largest absolute value of the entries in the constraint matrix $A$.
This assumption
%that will be made precise later, 
is based on the fact that the overwhelming majority of $[0,k]$-polytopes arising in combinatorial optimization for which an external description is known, satisfies $\alpha = 1$ \cite{SchBookCO}.
%
%More precisely, we design an algorithm for which the correspondin
%
%While only polynomially far from optimal, the bound 
%Our first contribution is the \emph{preprocessing \& scaling algorithm}.
%This 
%pivoting rule 
%simplex algorithm 
%generates a simplex path of length in $O(n^4 k \log (nk))$ for $[0,k]$-polytopes in $\R^n$.
%Based on this insight, we design a second simplex algorithm, named the \iterative.
Our second contribution is another simplex algorithm, named the \iterative, which generates a simplex path of length in $O(n^2k \log(nk\alpha))$. 
%The main difference with respect to the preprocessing \& scaling algorithm, is that we exploit the largest absolute value $\alpha$ of the entries in the constraint matrix $A$.
Thus, by exploiting the parameter $\alpha$, we are able to significantly improve the dependence on $n$.
%A crucial observation is that 
Moreover,
the dependence on $\alpha$ is only logarithmic, thus if $\alpha$ is bounded by a polynomial in $n, k$, then the length of our simplex path is in $O(n^2k \log(nk))$.  For $[0,1]$-polytopes this bound reduces to $O(n^2 \log n)$.
%
%\iffalse
%\note{$k=1$, $\alpha = 1$: $s-t$
%connector polytope,\\
%chain polytope, antichain polytope,\\
%matching polytope, perfect matching polytope, simple $b$-matching polytope, $b$-factor polytope\\
%vertex cover polytope, bipartite stable set polytope, independent set polytope\\
%edge cover polytope,\\
%transportation polytope,\\
%matchable set polytope,\\
%transversal polytope, common transversal polytope,\\
%$T$-join polytope\\
%forest polytope, branching polytope, arborescence polytope, matching forest polytope.\\
%Number of constraints exponential in the number of variables.
%\cite{SchBookCO}}
%\fi

%degrees of the polynomial that relates $\alpha$ to $n$ and $k$.
%$O(poly(n,k) \log(nk))$ since $K \in O(\alpha nk)$.
%Finally, for $[0,1]$-polytopes (where $k=1$) with $\alpha$ bounded by a polynomial in $n, k$, our bound becomes $O(n^2 \log(n))$, which improves over the bound of $O(n^3 k \log(nk))$ implied by Kitahara, Matsui and Mizuno.

%The second main contribution of this paper is the design of a
%The second key property of our pivoting rule is that
In both our simplex algorithms, our pivoting rule is such that
the number of operations needed to construct the next vertex in the simplex path is bounded by a polynomial in $n$, $m$, and $\log k$.
If $k$ is bounded by a polynomial in $n,m$, both our simplex algorithms are strongly polynomial.
This assumption is 
%motivated 
justified by the 
existence of $[0,k]$-polytopes that, for fixed $n$, have a diameter that grows almost linearly in $k$ \cite{DezPouSuk19}.
%fact that, as discussed in the previous paragraph, there exist $[0,k]$-polytopes that, for fixed $n$, have a diameter that grows almost linearly in $k$.
%following observations.
%As we have discussed in the previous paragraph, there are $[0,k]$-polytopes that, for fixed $n$, have a diameter that grows almost linearly in $k$.
Consequently, in order to obtain a simplex algorithm that is strongly polynomial also for these polytopes, we need to assume that $k$ is bounded by a polynomial in $n$ and $m$.
We remark that in this paper we use the standard notions regarding computational complexity in Discrete Optimization, and we refer the reader to Section~2.4 in the book~\cite{SchBookIP} for a thorough introduction.

\section{Overview of the algorithms and related work}

%Let $P$ be a $[0,k]$-polytope in $\R^n$.
%, i.e., a lattice polytope contained in $[0,k]^n$.

Our goal is to study the length of the simplex path
in the setting where the feasible region of problem \eqref{prob: main} is a $[0,k]$-polytope.
As discussed above, our main contribution is the design of a simplex algorithm that visits a number of vertices polynomial in $n$ and $k$, from any starting vertex $x^0$.
To the best of our knowledge, this question has not been previously addressed.

Since $P$ is a lattice polytope in $[0,k]^n$, it is not hard to see that a naive simplex algorithm, that we call the \basic, always reaches an optimal vertex of $P$ by constructing a simplex path of length at most $nk \norminf{c}$.
Therefore, in order to reach our goal, we need to make the simplex path length independent from $\norminf{c}$ .
%remove the dependence on $\norminf{c}$ in the length of the simplex path.

To improve the dependence on $\norminf{c}$ we design a \scaling that recursively calls the basic algorithm with finer and finer integral approximations of $c$. By using this bit scaling technique \cite{AhuMagOrl93}, we are able to construct a simplex path to an optimal vertex of $P$ of length in $O(nk {\log \norminf{c}})$.
 
The dependence of the simplex path length on $\norminf{c}$ is now logarithmic. Next, to completely remove the dependence on $\norminf{c}$, we first apply a \preprocessing due to Frank and Tardos \cite{FraTar87}, that replaces $c$ with a cost vector $\breve c$ such that
%$\log \norminf{\breve c}$ is polynomial in $n$ and $\log k$.
$\log \norminf{\breve c}$ is polynomially bounded in $n$ and $\log k$.
By using this preprocessing step in conjunction with our scaling algorithm, 
we design a \preprocessingscaling that constructs a simplex path whose length is polynomially bounded in $n$ and $k$.
\begin{theorem}
\label{th 1}
The \preprocessingscaling
generates a simplex path of length in\\$O(n^4k \log (nk))$.
\end{theorem}
%we are able to obtain an upper bound on the simplex path length that only depends on $n$ and $k$.

%Our second idea builds on the scaling algorithm with the goal of obtaining a simplex path of length polynomial in $n$ and $k$ alone, thus completely removing the dependence on $\norminf{c}$.
Our next task is to improve the polynomial dependence on $n$ by exploiting the knowledge of $\alpha$, the largest absolute value of an entry of $A$.
In particular, we design an \iterative which, at each iteration, identifies one constraint of $Ax \le b$ that is active at each optimal solution of \eqref{prob: main}.
Such constraint is then set to equality, effectively restricting the feasible region of \eqref{prob: main} to a lower dimensional face of $P$. 
At each iteration, we compute a suitable approximation $\tilde c$ of $c$ and we maximize $\tilde c ^{\top}x$ over the current face of $P$.
In order to solve this LP problem, we apply the \scaling and we trace a path along the edges of $P$.
We also compute an optimal solution to the dual, %of this LP problem, 
which is used to identify a new constraint of $Ax \le b$ that is active at each optimal solution of \eqref{prob: main}.
The final simplex path is then obtained by merging together the different paths constructed by the scaling algorithm at each iteration.

\begin{theorem}
\label{th 2}
%The length of the simplex path generated by the \iterative is in $O(n^2 k \log (nk\alpha))$.
The \iterative generates a simplex path of length in $O(n^2 k \log (nk\alpha))$.
%Given a  vertex $x^0 \in P$ and $c \in \Z^n$, the \iterative generates a simplex path $x^0, x^1,\dots, x^M$ in $P$ such that $x^M$ is an optimal solution to $\max\{c^{\top}x \mid x \in P\}$ and $M \in O(n^6k\log k)$. \cnote{it was $O(n^3k\log k)$}
\end{theorem}

%\cnote{output vector}

%At each iteration we compute a suitable approximation $\tilde c$ of $c$ and we apply the scaling algorithm to obtain a simplex path from the current vertex to one that maximizes $\tilde c^\top x$ over the current face of $P$. By computing an optimal solution to the dual of the latter LP, we are able to identify a new constraint of $Ax \le b$, active at the current vertex, that is active at each optimal solution of \eqref{prob: main}. The final simplex path is then obtained by merging together the different simplex paths constructed by the scaling algorithm at each iteration.

%\note{Mention that our simplex path is not improving wrt to $c$? I think it is now better to avoid this.}

%We remark that the simplex paths obtained with our scaling algorithm and with our iterative algorithm 

Our \iterative is inspired by Tardos' strongly polynomial algorithm for combinatorial problems \cite{Tar86}.
The three key differences with respect to Tardos' algorithm are:
%However, it differs from Tardos' algorithm in three aspects:
(1) Tardos' algorithm solves LP problems in standard form, while we consider a polytope $P$ in general form, i.e., $P = \{x \in \R^n \mid Ax \le b\}$; (2) Tardos' algorithm is parametrized with respect to the largest absolute value of a subdeterminant of the matrix A, while our algorithm relies on parameters $k$ and $\alpha$; and (3) Tardos' algorithm is \emph{not} a simplex algorithm, while our algorithm traces a simplex path along the edges of $P$.
Mizuno \cite{Miz16,MizSukDez18} proposed a variant of Tardos' algorithm that traces a \emph{dual} simplex path, under the assumption that $P$ is simple and that the constraint matrix is totally unimodular. We remark that the basic solutions generated by Mizuno's algorithm might not be primal feasible.

\smallskip

Next, we compare the length of the simplex path constructed by the \iterative to the upper bounds that are known for some other classic pivoting rules.
A result by Kitahara, Matsui and Mizuno \cite{KitMiz13,KitMatMiz12} implies that, if $Q = \{x \in \R^n_+ \mid Dx=d\}$ is a $[0,k]$-polytope in standard form with $D \in \Z^{p \times n}$ and $d \in \Z^p$, the length of the simplex path generated with Dantzig's or the best improvement pivoting rule is at most $(n-p) \cdot \min\{p,n-p\} \cdot k \cdot \log(k\min\{p,n-p\})$.
This has been recently extended to the steepest edge pivoting rule by Blanchard, De Loera and Louveaux \cite{BlaDeLLou20}.

Recall that in this paper we consider a lattice polytope $P \subseteq [0,k]^n$ in general form, i.e., $P = \{x \in \R^n \mid Ax \le b\}$. By applying an affine transformation, we can map $P$ to the standard form polytope $\bar P = \{(x,s) \in \R_+^{n + m} \mid Ax + I_m s = b\}$, where $I_m$ denotes the identity matrix of order $m$. Note that $\bar P$ is a lattice polytope and it has the same adjacency structure of $P$.
However, the affine transformation does not preserve the property of being a $[0,k]$-polytope.
%Thus, there is a one-to-one correspondence between the vertices of $P$ and those of $\bar P$, and two vertices are adjacent on $P$ is and only if their images are adjacent vertices of $\bar P$.
%vertex adjacencies are also preserved.
In fact, having introduced the {slack} variables $s$, we obtain that $\bar P$ is a lattice polytope in $[0,K]^{n+m}$, where $K = \max\{k, S\}$ and $S=\max\{\norminf{b-Ax} \mid x \in P\}$. In this setting, the result of Kitahara, Matsui and Mizuno implies that we can construct a simplex path whose length is
at most $n^2 K \log(nK)$,
%in $O(n^2 K \log(nK))$,
but $K$ critically depends on $S$, which in turns depends on $A,b$.
In fact, it is known that, even for $k=1$, the value $S$ can be as large as $\frac{(n-1)^{\frac{n-1}2}}{2^{2n+o(n)}}$ (see \cite{AloVu97,Zie00}), which is not polynomially bounded in $n,k$.
%Since, as mentioned in Section \ref{sec intro}, the diameter of a $[0,k]$-polytope in $\R^n$ is known to be at most $nk$,
%\cite{KleOnn92,dPMic16,DezPou18}, 
%this upper bound could be quite far from being tight.
%For $[0,k]$-polytopes, this upper bound could be quite far from being tight, as the diameter of $[0,k]$-polytopes is known to be at most $nk$ \cite{KleOnn92,dPMic16,DezPou18}. 
In turn, this implies that the length of the simplex path obtained in this way for $[0,k]$-polytopes is not polynomially bounded in $n,k$.

In particular, if each inequality of $Ax \le b$ is active at some vertex of $P$, we have that $S$ and $K$ are in $O( nk \alpha)$, thus in this case the upper bound implied by the result of Kitahara, Matsui and Mizuno is in $O(n^3 k \alpha \log( n k \alpha))$.
%Note that the dependence on $\alpha$ is superlinear, and even if $\alpha  = 1$, the upper bound is in $O(n^3 k \log( n k))$.
%In contrast, the dependence on $\alpha$ in the upper bound given by the \iterative is only logarithmic, and if $\alpha  = 1$ this upper bound is in $O(n^2k \log(nk))$.
%
Note that the dependence on $\alpha$ is superlinear, while in our \iterative it is only logarithmic.
%Consider now the setting $\alpha = 1$, which is the most favorable to 
Even if $\alpha  = 1$, the upper bound is in $O(n^3 k \log( n k))$. In contrast, the upper bound given by the \iterative is in $O(n^2k \log(nk))$.
% if $\alpha$ is bounded by a polynomial in $n, k$.
%
%is weaker than our bound of $O(n^2k \log(nk))$.
%Moreover, if $\alpha$ is bounded by a polynomial in $n, k$, then the length of the simplex path is in $O(n^2k \log(nk))$.
%
%This has important implications for $[0,1]$-polytopes with bounded $\alpha$. %Prominent examples of polytopes in this class are the matching polytope of nonbipartite graphs, the metric polytope for graphs with no $K_5$ minor and the stable set polytope of t-perfect graphs. Note that all these polytopes have an exponential number of constraints and are \emph{not} in standard form. In particular, bringing these $[0,1]$-polytopes into standard form would increase $k$ from $1$ to $O(n)$.\cnote{not correct. for matching and cut $k \in O(\sqrt n)$}.

To show that the bound of $O( nk \alpha)$ on $S$ and $K$ just discussed can be tight, we now provide an example of a $[0,1]$-polytope with $\alpha = 1$ and $S \in \Omega(n)$.
%As an example of a $[0,1]$-polytope with $\alpha = 1$
%\note{Can we take an odd cycle as the graph? Is it $t$-perfect? Otherwise (3) could be redundant}
Consider the stable set polytope of a $t$-perfect graph $G=(V,E)$, that is defined by the vectors $x \in \R^V_+$ satisfying:
\begin{align}
x_i + x_j &\le 1 && ij \in E \nonumber \\
\sum_{i \in V(C)}x_i &\le \floor{\frac{|V(C)|}{2}} && C \text{ odd cycle in $G$},\label{e: cycle inequalities}
%x \in \R^V_+
%x_i&\ge 0 & i\in V
\end{align}
where $V(C)$ denotes the nodes in the odd cycle $C$ \cite{SchBookCO}. Note that $x = 0$ is the characteristic vector of the empty stable set, thus it is a vertex of the stable set polytope. 
If $G$ is an odd cycle on $|V|=n$ nodes, then $G$ is $t$-perfect, and the inequality \eqref{e: cycle inequalities} corresponding to the cycle containing all nodes of $G$ is facet-defining.
Furthermore, the slack in such constraint can be as large as $\floor{\frac n 2}$, therefore $S \in \Omega(n)$. 
Consequently, the upper bound implied by \cite{KitMiz13,KitMatMiz12} is in $\Omega(n^3 \log n)$, while the upper bound given by our \iterative is in $O(n^2 \log n)$.

\section{The preprocessing \& scaling algorithm}
\label{s: th 1}

%\note{We remark that in this paper we use the standard notion of polynomial-time algorithms in Discrete Optimization, and we refer the reader to %Section~2.4 in
%the book~\cite{SchBookIP} for a thorough introduction.
%In particular, since we will be discussing strongly polynomial time algorithms, in this paper we use the arithmetic model of computation \cite{BerTsi97}.}

In the remainder of the paper, we study problem \eqref{prob: main}, where $c$ is a given cost vector in $\Z^n$, and $P$ is a $[0,k]$-polytope given via a system of linear inequalities, i.e., $P = \{x \in \R^n \mid Ax \le b\}$, where $A \in \Z^{m \times n}$, $b \in \Z^m$.
All our algorithms are simplex algorithms, meaning that they explicitly construct
%a path along the %reaches an optimal vertex by tracing 
a path along the edges of $P$ 
%of $P$ 
from \emph{any} given starting vertex $x^0$ to a vertex maximizing the liner function $c^\top x$.
%\cnote{Our goal is to explicitly construct a simplex path, i.e., a path along the edges of $P$ from the starting vertex $x^0$ to an optimal vertex $x^*$ of $c$. For this reason...}
%In this paper, we are interested in algorithms that are simplex-like, 
%We denote a path of this type as a \emph{simplex path.}
For this reason, we always assume that we are given a starting vertex $x^0$ of $P$.
%If one is interested in obtaining an arbitrary starting vertex $x^0$, it should be noted that it can be found in a number of operation polynomially bounded in $\size(A)$ via Tardos' algorithm.
Obtaining an arbitrary starting vertex $x^0$ can be accomplished via Tardos' algorithm by performing a number of operations that is polynomially bounded in $\size(A)$.
% \note{Justify this. Use  to find it. Write observation that size of A is polynomial in $n,m,\log k$.}
Recall that the \emph{size} of the matrix $A$, denoted by $\size(A)$, is in $O(nm \log \alpha)$ (see Section~2.1 in \cite{SchBookIP} for more details).

%In particular, the polytope is nonempty.
%Furthermore, we can assume without loss of generality that $P$ is full-dimensional \note{must remove this assumption} (see Lemma~1 in \cite{dPMic16}).
%Recall that, the polytope $P$ is 

The goal of this section is to present and analyze the \preprocessingscaling.
In order to do so, we first introduce the \oracle, that is the basic building block of all our algorithms, and then the \basic and the \scaling.
%The \scaling builds on the \basic.
%Recall that the input consists of a polytope $P = \{x \in \R^n \mid Ax \le b\}$, and a cost vector $c$.
%We recall that we always assume that the polyhedron $P$ is given via an external description, i.e., $P = \{x \in \R^n \mid Ax \le b\}$.

\subsection{The oracle}

%Clearly, also the polyhedron $P$ must be given in input. 
%To make our discussion more general, we do not assume that we have an explicit inequality description of $P$.

It will be convenient to consider the following \oracle, which provides a general way to construct the next vertex in the simplex path.
%\begin{itemize}
%\item[] 

%\medskip
%
%\noindent
%\textbf{Oracle.} \\
%\emph{Input:}
%A polytope $P$, a cost vector $c \in \Z^n$, and a vertex $x^t$ of $P$. 
%
%\noindent
%\emph{Output:}
%Either a statement that $x^t$ is optimal (i.e., $x^t \in \argmax \{{c}^\top x \mid x \in P\}$), or a vertex $x^{t+1}$ adjacent to $x^t$ with strictly larger cost (i.e., ${c}^\top x^{t+1} > {c}^\top x^t$). 
%
%
%\medskip

\begin{algorithm}[H]

\caption{\Oracle}
  
  \textbf{Input:} A polytope $P$, a cost vector $c \in \Z^n$, and a vertex $x^t$ of $P$.
  
  \textbf{Output:} Either a statement that $x^t$ is optimal (i.e., $x^t \in \argmax \{{c}^\top x \mid x \in P\}$), or a vertex $x^{t+1}$ adjacent to $x^t$ with strictly larger cost (i.e., ${c}^\top x^{t+1} > {c}^\top x^t$).

\end{algorithm}

%\note{replace $v$ with $x^t$ and $w$ with $x^{t+1}$?}

%
%\note{use $\ora$?}
%Instead, we assume that we have access to an \emph{oracle} that, if provided with a cost vector $c$ and a vertex $v$ of $P$, either returns a statement that $v$ is optimal (i.e., $v \in \argmin \{{c}^\top x \mid x \in P\}$), or returns a vertex $w$ adjacent to $v$ with strictly lower cost (i.e., ${c}^\top w < {c}^\top v$). 
%\note{Oracle complexity is $O(mn)$ only if nondegenerate}

In the next remark we analyze the complexity of the \oracle.
We recall that a polytope $P = \{x \in \R^n \mid Ax \le b\}$ is said to be \emph{simple} if, in each vertex of $P$, exactly $n$ inequalities from the system $Ax \le b$ are satisfied at equality.

\begin{remark}
\label{rem oracle complexity}
An \oracle call can be performed with a number of operations polynomially bounded in $\size(A)$. % \note{true or false?}.
If $P$ is simple, it can be performed in $O(nm)$ operations.
\end{remark}

\begin{prf}
If the polytope $P$ is simple, then the \oracle can be performed in $O(mn)$ operations.
This can be seen using a pivot of the dual simplex method, where the primal is in standard form, and the feasible region of the dual is given by the polytope $P$ \cite{BerTsi97}. 

Consider now the case where $P$ may be not simple.
Denote by $A^= x \le b^=$ the subsystem of the inequalities of $Ax \le b$ satisfied at equality by $x^t$.
Note that the polyhedron $T := \{x \in \R^n \mid A^= x \le b^=\}$ is a translated cone with vertex $x^t$.
Denote by $d^\top$ the sum of all the rows in $A^=$ and note that the vertex $x^t$ is the unique maximizer of $d^\top x$ over $T$. 
Let $T'$ be the truncated cone $T' := \{x \in T \mid d^\top x \ge d^\top x^t-1\}$ and note that there is a bijection between the neighbors of $x^t$ in $P$ and the vertices of $T'$ different from $x^t$.
We solve the LP problem $\max \{ c^\top x \mid x \in T'\}$.
Using Tardos' algorithm, this LP can be solved in a number of operations that is polynomial in the size of the constraint matrix, which is polynomial in $\size(A)$.

%Therefore, in order to 
%%obtain the neighbor $x^{t+1}$ of $x^t$ in $P$, 
%perform an oracle call,
%it suffices to solve the LP problem $\max \{ c^\top x \mid x \in T'\}$.
If $x^t$ is an optimal solution to the LP, then the \oracle returns that $x^t$ is optimal.
Otherwise, Tardos' algorithm returns an optimal solution that is a vertex $w$ of $T'$ different from $x^t$.
In this case the \oracle needs to return the corresponding neighbor $x^{t+1}$ of $x^t$ in $P$.
Let $A' x \le b'$ be the system obtained from $Ax \le b$ by setting to equality the inequalities in the subsystem $A^= x \le b^=$ satisfied at equality by both $x^t$ and $w$.
It should be noted that the vectors that satisfy $A' x \le b'$ constitute the edge of $P$ between $x^t$ and $x^{t+1}$.
The vector $x^{t+1}$ can then be found by maximizing $c^\top x$ over $A' x \le b'$ with Tardos' algorithm.
%The size of the constraint matrix is polynomial in $\size(A)$.
%
%
%In general case, keep only active constraints, construct cone, truncate it with objective that gives $x^t$, optimize our objective, find vertex. Bijection with neighbors of $x^t$.
%\note{A: finish this proof.}
\end{prf}

In all our algorithms, the simplex path is constructed via a number of \oracle calls with different inputs.
We note that, whenever $x^t$ is not optimal, our \oracle has the freedom to return \emph{any} adjacent vertex $x^{t+1}$ with strictly larger cost.
Therefore, our algorithms can all be customized by further requiring the \oracle to obey a specific pivoting rule.

%Our task in this paper is to obtain an algorithm to solve \eqref{prob: main} %which performs as few oracle calls as possible.
%via a number of oracle calls.
%The number of oracle calls coincides with the length of the generated simplex path.
%The algorithms to solve \eqref{prob: main} that we introduce in this paper will all invoke this oracle as a subroutine.
%In this paper we  is to obtain an algorithm  %which performs as few oracle calls as possible.
%via a number of oracle calls.

\subsection{The basic algorithm}

The simplest way to solve \eqref{prob: main} is to
%Note that we could simply recursively  
recursively invoke the \oracle with input $P$, $c$, and the vertex $x^t$ obtained from the previous iteration, starting from the vertex $x^0$ in input.
We formally describe this \basic, which will be used as a subroutine in our 
%\scaling described later in the section.
%subsequent algorithms.
%more complex algorithms that will be introduced later.
subsequent algorithms.

%\medskip
%
%\noindent
%\textbf{Basic algorithm.} \\
%\emph{Input:} 
%A $[0,k]$-polytope $P$, a cost vector $c \in \Z^n$, and a vertex $x^0$ of $P$.
%
%\noindent
%\emph{Output:} A vertex of $P$ maximizing $c^\top x$. \\
%For $t =0,1,2,\dots$: 
%Invoke the oracle with input $P$, $c$, and $x^t$.
%If the output of the oracle is a statement that $x^t$ is optimal, return $x^t$.
%Otherwise, let $x^{t+1}$ be the vertex of $P$ returned by the oracle.
%
%\medskip

\begin{algorithm}[H]

\caption{\Basic}

\textbf{Input:} A $[0,k]$-polytope $P$, a cost vector $c \in \Z^n$, and a vertex $x^0$ of $P$.

\textbf{Output:} A vertex $x^*$ of $P$ maximizing $c^\top x$.

\begin{algorithmic}
\For{$t =0,1,2,\dots$}
\State Invoke $\oracle(P, c, x^t)$.
\State If the output of the \oracle is a statement that $x^t$ is optimal, return $x^t$.
\State Otherwise, set $x^{t+1} := \oracle(P, c, x^t)$.
\EndFor
\end{algorithmic}
\end{algorithm}

The correctness of the \basic is immediate.
Next, we upper bound the
%the number of oracle calls performed by the \basic.
length of the simplex path generated by the \basic.
%In the remainder of this section we denote by $\gamma := \norminf{c} = \max\{\abs{c_1}, \dots, \abs{c_n}\}$.

%\note{replace $U$ with greek letter $\gamma$ (like we have $\alpha$ for $A$) for consistency.}

\begin{observation}
\label{obs basic}
%The \basic calls the oracle at most $nk\norminf{c}$ times.
The length of the simplex path generated by the \basic is bounded by 
$c^\top x^* - c^\top x^0$. %, where $x^*$ is the vertex of $P$ in output.
In particular, it is bounded by $nk\norminf{c}$.
\end{observation}

\begin{prf}
%This bound follows from the observation that 
To show the first part of the statement, we only need to observe that each \oracle call increases the objective value by at least one, since $c$ and the vertices of $P$ are integral.

%To prove the second part of the statement, we show that 
The cost difference % $c^\top x^* - c^\top x^0$ 
between $x^*$ and $x^0$ of $P$ can be bounded by 
\[
c^\top x^* - c^\top x^0 = 
\sum_{i=1}^n c_i (x^*_i-x^0_i)
\le \sum_{i=1}^n \abs{c_i} \abs{x^*_i-x^0_i} 
\le  nk \norminf{c},
\]
where the last inequality we use $\abs{x^*_i-x^0_i} \le k$ since $P$ is a $[0,k]$-polytope.
This concludes the proof of the second part of the statement.
%We only need to show that the \basic invokes the oracle at most $nk\norminf{c}$ times.
\end{prf}
%\note{In the above lemma we could have used the L1 norm instead to obtain a better bound. Correct? Should we? This alternative bound should be $k V$, where $V := \normone{c} = \sum_{i=1}^n \abs{c_i}$.}

%Observation \ref{obs basic} implies that the path generated by the \basic has length at most $nk\norminf{c}.$

From Remark~\ref{rem oracle complexity}, we directly obtain the following.

%\note{CHECK}

\begin{remark}
\label{rem basic complexity}
The number of operations performed by the \basic to construct the next vertex in the simplex path is polynomially bounded in $\size(A)$.
If $P$ is simple, the number of operations is $O(nm)$.
\end{remark}

\subsection{The scaling algorithm}

%In particular, we do not seem to be able to give any bound on the number of \oracle calls performed by the \basic that is polynomial in $k,n,\log \norminf{c}$.
%\note{Sure? What about using the Onn construction? Do we obtain an upper bound on number of vertices? It would be great if we could give example that $nk\norminf{c}$ is tight.}

The length of the simplex path generated by the \basic is clearly not satisfactory.
%This is in stark contrast with the fact that 
In fact, as we discussed in Section~\ref{sec intro}, our goal is to obtain a simplex path of length polynomial in $n$ and $k$, and therefore 
%the diameter of $P$ is bounded by $O(nk)$, which in particular is 
independent from $\norminf{c}$.
In this section we improve this gap by giving a \scaling that 
%solves problem \eqref{prob: main} by invoking the \oracle at most $nk (\ceil{\log \norminf{c}} + 1)$ times.
%In particular, this algorithm 
yields a simplex path 
%from the input vertex $x^0$ to an optimal solution of \eqref{prob: main} 
of length in $O(nk \log \norminf{c})$.
%We will refer to this algorithm as the \emph{\scaling}.

Our \scaling is based on a \emph{bit scaling} technique.
For ease of notation, we define $\ell := \ceil{\log \norminf{c}}$.
The main idea is to iteratively use the \basic with the sequence of increasingly accurate integral approximations of the cost vector $c$ given by
%For $t=0,\dots, \ell$, let $c^t$ be the integral vector defined as
\[
c^t := \ceil{\frac{c}{2^{\ell-t}}} \qquad \text{for $t=0,\dots, \ell$.}
\]
Since $c$ is an integral vector, we have $c^\ell = c$.

Bit scaling techniques have been extensively used to develop polynomial-time algorithms for a wide array of discrete optimization problems.
Edmonds and Karp \cite{EdmKar72} and Dinic \cite{Din73} independently introduced this technique in the context of the minimum cost flow problem.
%known in the field of optimization since at least the 1970s.
Gabow \cite{Gab85} used it for shortest path, maximum flow, assignment, and matching problems, 
while Schulz, Weismantel, and Ziegler \cite{SchWeiZie95} 
employed it to design a primal algorithm for $0/1$-integer programming.
%designed a primal algorithm for $0/1$-integer programming based on scaling.
The book \cite{AhuMagOrl93} popularized bit scaling as a generic algorithmic tool in optimization.

%minimum cost flow \cite{Roc80}

%We refer the reader to Chapter~3 in \cite{AhuMagOrl93} for more 

%Scaling methods have extensivelybeen used to derive polynomial time al- gorithms for a wide variety of network and combinatorial optimization prob- lems (see, e.g., [AMO93]). In this section we use bit-scaling of costs to derive an \oracle-polynomialtime algorithm for optimizing 0/l-integer programs.This techniquehas beenusedearlier by RSck [R580] and Gabow [Gab85] for solving minimum cost flow and shortestpath problems,respectively.

Next, we describe our algorithm.

%\medskip
%
%\noindent
%\textbf{Scaling algorithm.} \\
%\emph{Input:} A $[0,k]$-polytope $P$, a vertex $x^0$ of $P$, and a cost vector $c \in \Z^n$. 
%
%
%\noindent \emph{Output:} A vertex of $P$ maximizing $c^\top x$. \\
%For $t = 0,\dots,\ell$:
%Invoke the \basic 
%with input $P$, $c^t$, and $x^t$.
%Let $x^{t+1}$ be the vertex of $P$ returned by the \basic. \\
%Return the vertex $x^{\ell+1}$.
%
%
%\medskip

\begin{algorithm}[H]

\caption{\Scaling}

  \textbf{Input:} A $[0,k]$-polytope $P$, a cost vector $c \in \Z^n$, and a vertex $x^0$ of $P$.
  
  \textbf{Output:} A vertex $x^*$ of $P$ maximizing $c^\top x$.

\begin{algorithmic}
\For{$t = 0,\dots,\ell$}
%\State Call the \basic with input $P$, $c^t$, and $x^t$.
%\State Let $x^{t+1}$ be the vertex of $P$ returned by the \basic.
\State Compute $c^t$.
\State Set $x^{t+1} := \basic(P, c^t, x^t)$.
\EndFor

\State Return $x^{\ell+1}$.
\end{algorithmic}
\end{algorithm}

The correctness of the \scaling follows from the correctness of the \basic, since the vector $x^{\ell+1}$ returned is the output of the \basic with input $P$ and cost vector $c^\ell = c$.
%Hence, the vector $x^{\ell+1}$ is a vertex of $P$ maximizing ${c^\ell}^\top x = c^\top x$.

Next, we analyze the length of the simplex path generated by the \scaling.
We first derive some simple properties of the approximations $c^t$ of $c$.

\begin{lemma}
\label{lem: boxi}
For each $t=0,\dots, \ell$, we have $\norminf{c^t} \le 2^t$.
\end{lemma}

\begin{prf}
By definition of $\ell$, we have $\abs{c_j} 
\le \norminf{c} 
\le 2^\ell$ for every $j=1,\dots,n$, hence $-2^\ell \le c_j \le 2^\ell$.
For any $t \in \{0,\dots, \ell\}$, we divide the latter chain of inequalities by $2^{\ell-t}$ and round up to obtain 
\[
- 2^t = \ceil{- 2^t} 
= \ceil{\frac{-2^\ell}{2^{\ell-t}}} 
\le \ceil{\frac{c_j}{2^{\ell-t}}} 
\le \ceil{\frac{2^\ell}{2^{\ell-t}}} 
= \ceil{2^t} = 2^t.
\]
\end{prf}
\begin{lemma}
\label{lem: dist}
For each $t = 1,\dots, \ell$, we have $2c^{t-1} - c^t \in \{0, 1\}^n$.
\end{lemma}
\begin{prf}
First, we show that for every real number $r$, we have $2 \ceil{r} - \ceil{2 r} \in \{0, 1\}$.
Note that $r$ can be written as $\ceil{r} + f$ with $f \in (-1,0]$.
We then have $\ceil{2 r} = \ceil{2 \ceil{r} + 2 f} = 2 \ceil{r} + \ceil{2 f}$.
Since $\ceil{2 f} \in \{-1,0\}$, we obtain $\ceil{2 r} - 2 \ceil{r} \in \{-1, 0\}$, hence $2 \ceil{r} - \ceil{2 r} \in \{0, 1\}$.

Now, let $j \in \{1,\dots,n\}$, and consider the $j$th component of the vector $2c^{t-1} - c^t$.
By definition, we have 
\[
2c^{t-1}_j - c^t_j 
= 2 \ceil{\frac{c_j}{2^{\ell-t+1}}} - \ceil{\frac{c_j}{2^{\ell-t}}}.
\]
The statement then follows from the first part of the proof by setting $r = c_j / 2^{\ell-t+1}$.
\end{prf}

We are ready to provide our bound on the length of the simplex path generated by the \scaling.
Even though the \scaling uses the \basic as a subroutine, we show that the simplex path generated by the \scaling is much shorter than the one generated by the \basic alone.

\begin{proposition}
\label{prop scaling oracle calls}
The length of the simplex path generated by the \scaling is bounded by $nk (\ceil{\log \norminf{c}} + 1) \in O(nk \log \norminf{c})$.
\end{proposition}

\begin{prf}
%We only need to show that the \scaling invokes the \oracle at most $nk (\ceil{\log \norminf{c}} + 1)$ times.
Note that the \scaling performs a total number of $\ell+1 = \ceil{\log \norminf{c}} + 1$ iterations, and in each iteration it calls once the \basic.
Thus, we only need to show that, at each iteration, the simplex path generated by the \basic is bounded by $nk$.

First we consider the iteration $t =0$ of the \scaling.
In this iteration, the \basic is called with input $P$, $c^0$, and $x^0$.
%used to solve the linear program $\min \{{c^0}^\top x \mid x \in P\},$ starting from the input vertex $x^0$.
Lemma~\ref{lem: boxi} implies that $\norminf{c^0} \le 1$, and from Observation~\ref{obs basic} we have that the \basic calls the \oracle at most $nk$ times.

Next, consider the iteration $t$ of the \scaling for $t \in \{1,\dots,\ell\}$.
In this iteration, the \basic is 
%used to solve the linear program \eqref{prob: aux} starting from the vertex $x^t \in \argmin \{{c^{t-1}}^\top x \mid x \in P\}$.
called with input $P$, $c^t$, and $x^t$, and outputs the vertex $x^{t+1}$.
From Observation \ref{obs basic}, we only need to show that ${c^t}^\top x^{t+1} - {c^t}^\top x^t \le nk.$

First, we derive an upper bound on ${c^t}^\top x^{t+1}$.
By construction of $x^t$, the inequality ${c^{t-1}}^\top x \le {c^{t-1}}^\top x^t$ is valid for the polytope $P$, thus
\[
{c^t}^\top x^{t+1} = \max \{ {c^t}^\top x \mid x \in P \} \le \max \{ {c^t}^\top x \mid x \in [0,k]^n, \ {c^{t-1}}^\top x \le {c^{t-1}}^\top x^t \}. 
\]
The 
%maximization problem on the 
optimal value of the LP problem on the 
right-hand side is upper bounded by $2{c^{t-1}}^\top x^t$.
In fact, Lemma~\ref{lem: dist} implies $c^{t} \le 2 c^{t-1}$, hence for every feasible vector $x$ of the LP problem on the right-hand side, we have
\[
{c^{t}}^\top x \le 2{c^{t-1}}^\top x \le 2{c^{t-1}}^\top x^t.
\]
Thus we have shown ${c^t}^\top x^{t+1} \le 2{c^{t-1}}^\top x^t$.

We can now show ${c^t}^\top x^{t+1} - {c^t}^\top x^t \le nk.$
%We can now upper bound the quantity ${c^t}^\top x^{t+1} - {c^t}^\top x^t$.
We have
\[
{c^t}^\top x^{t+1} - {c^t}^\top x^t 
\le 2{c^{t-1}}^\top x^t - {c^t}^\top x^t
= (2c^{t-1} - c^{t})^\top x^t \le nk.
\]
The last inequality holds because, from Lemma~\ref{lem: dist}, we know that $2c^{t-1} - c^{t} \in \{0, 1\}^n$, while the vector $x^t$ is in $[0,k]^n$.
%This concludes the proof of the proposition.
%At iteration $t$ of the \scaling, the \basic solves the problem \eqref{prob: aux} starting from the vertex $x^t$ and eventually reaching the vertex $x^{t+1}$.
%Observation \ref{obs basic} then implies that the length of the simplex path generated at iteration $t$ is bounded by $nk$.
%Since each \oracle call decreases the objective value by at least one, and since the vector $c^t$ and the vertices of $P$ are integral, we obtain that the \basic calls the \oracle at most $nk$ times.
\end{prf}

\begin{remark}
\label{rem scaling complexity}
The number of operations performed by the \scaling 
to construct the next vertex in the simplex path
is polynomially bounded in $\size(A)$ and $\log \norminf{c}$.
If $P$ is simple, the number of operations is in $O(n^2m \log^2 \norminf{c})$.
\end{remark}

\begin{prf}
Let $x^i$ be the $i$-th vertex of the simplex path computed by the \scaling.
%First, consider the case where $x^i$ is not the optimal vertex returned by the \basic in some iteration of the \scaling.
%In this case, the next vertex $x^{i+1}$ of the simplex path is obtained at the next iteration of the \basic, which requires only one \oracle call.
%
%Next, consider the case where $x^i$ is the optimal vertex returned by the \basic in some iteration of the \scaling.
%The algorithm then restarts the for cycle, it computes an approximation $c^t$ of $c$,
%and then calls once again the \basic.
%
%At this point, we distinguish two sub-cases.
%In the first sub-case, $x^i$ is not optimal for the new LP with cost $c^t$.
%Then, the first \oracle call of the \basic yields the next vertex $x^{i+1}$ in the simplex path. 
The \scaling might call, in the worst case, $\ell$ times the \basic before obtaining the next vertex $x^{i+1}$.
Each time, the \scaling first computes an approximation $c^t$ of $c$ and then calls the \basic. % which, in turn, invokes the \oracle only once.
Computing $c^t$ can be done by binary search, and the number of comparisons required is at most $n\log \norminf{c^t}$, which is bounded by $nt \le n \ell$ from Lemma~\ref{lem: boxi}.
Furthermore, from Remark~\ref{rem basic complexity}, each time the \basic is called, it performs a number of operations polynomially bounded in $\size(A)$, and by $O(nm)$ if $P$ is simple.
Therefore, to compute $x^{i+1}$ we need 
a number of operations bounded a polynomial in $\size(A)$ and in $\log \norminf{c}$.
If $P$ is simple, the number of operations is $O(\ell \cdot n\ell \cdot nm) = O(n^2m \log^2 \norminf{c})$.
\end{prf}

In the next section we use the \scaling as a subroutine in the \preprocessingscaling.
We remark that, the \scaling will also be a subroutine in the \iterative, which is described in Section~\ref{sec th 2}.
% will be used as a subroutine in our two main algorithms, which are the \preprocessingscaling, given in Section~\ref{sec preprocessing and scaling}, and the 

%%%%%%%%%%%%%%%%%%%%%%%%%%%%%%%%%%%%%%%%%%%%%%%%%%%%%
\subsection{The preprocessing \& scaling algorithm}
\label{sec preprocessing and scaling}

The length of the simplex path generated by the \scaling still depends on $\norminf{c}$, even though the dependence is now logarithmic instead of linear.
%This is in stark contrast with the fact that 
%In this section we present our \iterative, which completely removes the dependence on $\norminf{c}$.
In this section we show that we can completely remove the dependence on $\norminf{c}$ by using  our \scaling in conjunction with the \preprocessing by 
%We remark that the \scaling will be used as a subroutine in the \iterative.
%In fact, as we discussed in Section~\ref{sec intro}, our goal is to obtain a simplex path of length polynomial in $n$ and $k$, and therefore 
%the diameter of $P$ is bounded by $O(nk)$, which in particular is 
%independent on $\norminf{c}$.
%
%
%
Frank and Tardos \cite{FraTar87}.
This method relies on the simultaneous approximation algorithm of Lenstra, Lenstra and Lov\'asz \cite{LenLenLov82}.
%Quite recently Fujishige (1985) developed a version of the algorithm that also applies to situations in which A has exponentially many rows.
Next, we state the input and output of Frank and Tardos' algorithm.

%\medskip
%
%\noindent
%\textbf{Preprocessing algorithm.} \\
%\emph{Input:} A cost vector $c \in \Q^n$, and $N \in \Z$.
%
%
%\noindent \emph{Output:} A cost vector $\breve c \in \Z^n$ such that $\norminf{\breve c} \le 2^{4n^3} N^{n(n+2)}$ and $\sign(c z) = \sign(\breve c z)$ for every $z \in \Z^n$ with $\normone{z} \le N-1$. \\
%
%\medskip

\begin{algorithm}[H]

\caption{\Preprocessing}

  \textbf{Input:} A vector $c \in \Q^n$ and a positive integer $N$.

  \textbf{Output:} A vector $\breve c \in \Z^n$ such that $\norminf{\breve c} \le 2^{4n^3} N^{n(n+2)}$ and $\sign(c^\top z) = \sign(\breve c^\top z)$ for every $z \in \Z^n$ with $\normone{z} \le N-1$.

\end{algorithm}

The number of operations performed by the \preprocessing is polynomially bounded in $n$ and $\log N$.
For more details, we refer the reader to Section 3 in \cite{FraTar87}.

%\begin{theorem}[Theorem 3.3 in \cite{FraTar87}]
%Let $c \in \Q^n$, and $N \in \Z$.
%There is an algorithm that computes a vector $\breve c \in \Z^n$ such that $\norminf{\breve c} \le 2^{4n^3} N^{n(n+2)}$ and $\sign(c z) = \sign(\breve c z)$ for every $z \in \Z^n$ with $\normone{z} \le N-1$.
%The number of operations performed by this algorithm is polynomially bounded in $n$ and $\log N$.
%\end{theorem}

Next, we describe the algorithm obtained by combining the \preprocessing and the \scaling.

%\medskip
%
%\noindent
%\textbf{Preprocessing \& scaling algorithm.} \\
%\emph{Input:} A $[0,k]$-polytope $P$, a vertex $x^0$ of $P$, and a cost vector $c \in \Z^n$. 
%
%
%\noindent \emph{Output:} A vertex of $P$ maximizing $c^\top x$. \\
%Invoke the preprocessing algorithm with input $c$ and $N := 2k$, and let $\breve c$ be the vector returned. \\
%Invoke the scaling algorithm with input $P$, $\breve c$, and $x^0$, and let $x^1$ be the returned vertex of $P$. \\
%Return $x^1$.
%
%\medskip

\begin{algorithm}[H]

\caption{\Preprocessingscaling}
  
  \textbf{Input:} A $[0,k]$-polytope $P$, a cost vector $c \in \Z^n$, and a vertex $x^0$ of $P$.
  
  \textbf{Output:} A vertex $x^*$ of $P$ maximizing $c^\top x$.

\begin{algorithmic}
\State Set $\breve c := \preprocessing(c,N := nk+1)$.
\State Set $x^* := \scaling(P, \breve c, x^0)$.
\State Return $x^*$.
\end{algorithmic}
\end{algorithm}

We first show that the \preprocessingscaling is correct.

\begin{proposition}
\label{prop preprocessing and scaling correct}
The vertex $x^*$ of $P$ returned by the \preprocessingscaling maximizes $c^\top x$ over $P$.
\end{proposition}

\begin{prf}
Due to the correctness of the \scaling, we have that $\breve c^\top (x^*- x) \ge 0$ for every $x \in P$.
Note that, for every $x \in P \cap \Z^n$, we have $x^*- x \in \Z^n$ and $\normone{x^*- x} \le nk = N-1$.
Therefore, the \preprocessing guarantees that $c^\top (x^*- x) \ge 0$ for every $x \in P \cap \Z^n$.
The statement follows because all vertices of $P$ are integral.
\end{prf}

We are now ready to give a proof of Theorem~\ref{th 1}.
We show that the obtained simplex path length is polynomially bounded in $n,k$, thus only polynomially far from the worst-case diameter.

%\begin{proposition}
%\label{prop preprocessing and scaling oracle calls}
%The length of the simplex path generated by the \preprocessingscaling is bounded by $nk (4n^3 \log 2 + n(n+2) \log (nk+1) + 2) \in O(n^4 k \log (nk))$.
%\end{proposition}

\begin{prfc}[of Theorem~\ref{th 1}]
The vector $\breve c$ returned by the \preprocessing satisfies $\norminf{\breve c} \le 2^{4n^3} (nk+1)^{n(n+2)}$, hence $\log \norminf{\breve c} \le 4n^3 \log 2 + n(n+2) \log (nk+1)$.
From Proposition \ref{prop scaling oracle calls}, the length of the simplex path generated by the \preprocessingscaling is bounded by 
\begin{align*}
nk (\ceil{\log \norminf{\breve c}} + 1) 
\le nk (4n^3 \log 2 + n(n+2) \log (nk+1) + 2)  \in O(n^4 k \log (nk)).
\end{align*}
\end{prfc}

We conclude this section by analyzing the number of operations performed by the \preprocessingscaling.

\begin{remark}
\label{rem preprocessing scaling complexity}
The number of operations performed by the \preprocessingscaling 
to construct the next vertex in the simplex path
is polynomially bounded in $\size(A)$ and $\log k$.
If $P$ is simple, the number of operations is polynomially bounded in $n,m,$ and $\log k$.
\end{remark}

\begin{prf}
The number of operations performed by the \preprocessingscaling to construct the next vertex in the simplex path is the sum of: 
(i) the number of operations needed to compute $\breve c$, and
(ii) the number of operations performed by the \scaling, with cost vector $\breve c$, to construct the next vertex in the simplex path.
The vector $\breve c$ can be computed with a number of operations polynomially bounded in $n$ and $\log(nk)$ \cite{FraTar87}.
From Remark~\ref{rem scaling complexity}, (ii) is polynomially bounded in $\size(A)$ and $\log \norminf{\breve c}$, and by $O(n^2m \log^2 \norminf{\breve c})$ if $P$ is simple.
To conclude the proof, we only need to observe that $\log \norminf{\breve c}$ is polynomially bounded in $n$ and $\log k$.
\end{prf}

%Therefore, the \preprocessingscaling yields a simplex path of length $O(n^4 k + kn^3 \log (nk))$.

%%%%%%%%%%%%%%%%%%%%%%%%%%%%%%%%%%%%%%%%%%%%%%%%%%%%%
\section{The iterative algorithm}
\label{sec th 2}

%\note{Notation problem: The letter P is used for primal and the polytope. Can we just erase Primal and Dual here? They are never used in the algorithm.}

%We further assume that each inequality in $Ax \le b$ is facet-defining \note{still needed?}, and that 
%%, without loss of generality, that 
%the greatest common divisor of the entries in each row of $A$ is one \note{still needed?}.
%Both these assumptions are without loss of generality and it is well-known that we can reduce ourselves to this setting in polynomial time.
%%inequality $a_i^\top x \le b_i$ in $Ax \le b$
In this section, we design a simplex algorithm that yields simplex path whose length depends on $n,k$, and $\alpha$, where $\alpha$ denotes the largest absolute value of the entries of $A$.
We define $[m] := \{1,2,\dots,m\}$ and refer to the rows of $A$ as $a_i$, for $i \in [m]$.
Next, we present our \iterative.

\begin{algorithm}[H]
  
\caption{\Iterative}
  
\textbf{Input:} A $[0,k]$-polytope $P$, a cost vector $c \in \Z^n$, and a vertex $x^0$ of $P$.
  
\textbf{Output:} A vertex $x^*$ of $P$ maximizing $c^\top x$.

\begin{algorithmic}[1]

\setcounter{ALG@line}{-1}

\State 
Let $\E := \emptyset$ and $x^* := x^0$.
\label{alg: init}

\State 
Let $\bar c$ be the projection of $c$ onto the subspace $\{x\in \R^n  \mid a_i^\top x = 0 \text{ for } i \in \E\}$ of $\R^n$. 
If $\bar c = 0$ return $x^*$, otherwise go to~\ref{alg: 2}.
\label{alg: 1}

\State 
Let $\tilde c \in \Z^n$ be defined by $\tilde c_i := \floor{\frac{n^3 k \alpha}{\norminf{\bar c}} \bar c_i}$ for $i = 1,\dots, n$.
\label{alg: 2}

\State 
Consider the following pair of primal and dual LP problems:
\label{alg: 3}

\noindent\begin{minipage}{0.5\linewidth}
\begin{equation}
\begin{array}{rll}
\max & \tilde c^\top x&\\
\text{s.t.~} & a_i^\top x = b_i & i \in \E \\
& a_i^\top x \le b_i & i \in [m]\setminus \E
\end{array}\tag{$\tilde P$}\label{tilde P}
\end{equation}
\end{minipage}%
\begin{minipage}{0.5\linewidth}
\begin{equation}
\begin{array}{rll}
\min & b^\top y& \\
\text{s.t.~} & A^\top y  = \tilde c&\\
& y_i \ge 0 & i \in [m]\setminus \E.
\end{array}\tag{$\tilde D$}\label{tilde D}
\end{equation}
\end{minipage}\par\vspace{\belowdisplayskip}

\noindent
Use the \scaling to compute an optimal vertex $\tilde x$ of \eqref{tilde P} starting from $x^*$.

\noindent
Compute an optimal solution $\tilde y$ to the dual \eqref{tilde D}
such that
(i)
$\tilde y$ has at most $n$ nonzero components, and
(ii)
$\tilde y_j = 0$ for every $j \in [m] \setminus \E$ such that $a_j$ can be written as a linear combination of $a_i$, $i \in \E$.

\noindent
Let $\H := \{i \mid \tilde y_i > n k \}$, and let $h \in \H \setminus \E$.
Add the index $h$ to the set $\E$, set $x^* := \tilde x$, and go back to step~\ref{alg: 1}.

\end{algorithmic}
\end{algorithm}

%\cnote{I would consider modifying how step 3 is written to be more in the spirit of the simplex method. Example:
%Let $\E(\tilde x) \supset \E$ be the indices of the constraints of \eqref{tilde P} active at $\tilde x$. Let $\B \subset \E(\tilde x)$, be such that: (i)$\abs{\B}=n$ and $ \E \subset \B $, (ii) the row submatrix $A_{\B}$ of $A$ indexed by $\B$ is nonsingular and (iii) the $n$-dimensional vector $y_{\B} = (A_{\B}^{\top})^{-1} \tilde c$ whose entries are indexed by $\B$ is such that $y_i \ge 0$ for $i \in \E(\tilde x)$}.

As the name suggests, the above algorithm is iterative in nature.
In particular, an \emph{iteration} of the algorithm corresponds to one execution of steps~\ref{alg: 1},~\ref{alg: 2}, and~\ref{alg: 3}.

\subsection{Well-defined}

In this section we prove that the \iterative is well-defined, meaning that it can indeed perform all the instructions stated in its steps. % problem \eqref{tilde P} is always feasible.
First we show that, when the \iterative calls the \scaling in step~\ref{alg: 3}, $x^*$ is a valid input.
%First, we show that the \iterative can indeed call the \scaling in step~\ref{alg: 3}.

\begin{lemma}
\label{l: feasible}
%The \iterative can call the \scaling in step~\ref{alg: 3} with input vertex $x^*$.
In step~\ref{alg: 3}, the \iterative calls the \scaling with valid inputs. 
\end{lemma}

\begin{prf}
We only need to show that the vector $x^*$ is feasible for problem \eqref{tilde P}.
Consider the vectors $\tilde x, \tilde y$ and the index $h$ from the previous iteration.
We have $\tilde y_h > 0$, which, by complementary slackness, implies $a_h^\top \tilde x = b_h$.
%To see this, one just needs to observe that, in step~\ref{alg: 3}, the vector $x^*$ 
Therefore, the $\tilde x$ from the previous iteration, which coincides with $x^*$ of the current iteration, is indeed feasible for the problem \eqref{tilde P} of the current iteration.
\end{prf}

%\note{Do we ever show that we obtain a simplex path? This is because $\tilde x$ is in the face of the next iteration, because $\tilde y_i >0$ and complementary slackness.}

Next, we show that an optimal solution $\tilde y$ to the dual \eqref{tilde D} with the properties stated in step~\ref{alg: 3} exists, and it can be computed efficiently.

\begin{lemma}
\label{l: compute dual vars}
In step~\ref{alg: 3} of the \iterative, a vector $\tilde y$ satisfying (i) and (ii) always exists, and the number of operations needed to compute it is polynomially bounded in $\size(A)$. 
If $P$ is simple, the number of operations is in $O(nm+n^3)$.
\end{lemma}

\begin{prf}
First, assume that the polytope $P$ is simple.
Let problem \eqref{tilde P}' be obtained from \eqref{tilde P} by dropping the inequalities that are not active at $\tilde x$.
This can be done in $O(nm)$ operations.
Since $P$ is simple, the number of constraints in \eqref{tilde P}' is $n$ and the $n \times n$ constraint matrix $\tilde A$ of \eqref{tilde P}' is invertible.
 Note that $\tilde x$ is an optimal solution to \eqref{tilde P}' as well.
 Let \eqref{tilde D}' be the dual of \eqref{tilde P}'.
 Note that \eqref{tilde D}' is obtained from \eqref{tilde D} by dropping the variables $y_j$ corresponding to the inequalites of \eqref{tilde P} dropped to obtain \eqref{tilde P}'.
 Since \eqref{tilde P}' has an optimal solution, then so does \eqref{tilde D}' from strong duality.
 The constraint matrix of \eqref{tilde D}' is the invertible matrix $\tilde A^\top$.
 The only feasible solution to the system of linear equations in \eqref{tilde D}' is the vector $\tilde y' := \tilde A^{-\top} \tilde c$ which can be computed in $O(n^3)$ operations.
 Since \eqref{tilde D}' is feasible, then $\tilde y'$ must satisfies all the constraints in \eqref{tilde D}', thus $\tilde y'$ is optimal.
 Let $\tilde y$ be obtained from $\tilde y'$ by adding back the dropped components and setting them to zero.
 The vector $\tilde y$ is feasible to \eqref{tilde D}, and, from complementary slackness with $\tilde x$, it is optimal to \eqref{tilde D}.
 Furthermore, $\tilde y$ clearly satisfies (i).
 To see that it satisfies (ii), note that the inequalities 
 $a_i^\top x = b_i$, $i \in \E$, are all in \eqref{tilde P}'.
 Since the constraints in \eqref{tilde P}' are all linearly independent, problem \eqref{tilde P}' cannot contain any constraint $a_j^\top x \le b_j$, for $j \in [m] \setminus \E$ such that $a_j$ can be written as a linear combination of $a_i$, $i \in \E$.
 Hence, the corresponding dual variable $\tilde y_j$ has been set to zero.

 \smallskip

Consider now the general case where $P$ may be not simple.
First, we show how to compute a vector $\tilde y$ that satisfies (i).
Since \eqref{tilde P} has an optimal solution, then so does \eqref{tilde D} from strong duality.
%Note that \eqref{tilde D} is equivalent to 
%%\noindent\begin{minipage}{0.5\linewidth}
%%\begin{equation}
%%\begin{array}{rll}
%%\max & \tilde c^\top x&\\
%%\text{s.t.~} & a_i^\top x \le b_i & i \in [m]\setminus \E\\
%%& a_i^\top x = b_i & i \in \E.
%%\end{array}\tag{$\textcolor{red}{\tilde P}$}\label{tilde P}
%%\end{equation}
%%\end{minipage}%
%%\begin{minipage}{0.5\linewidth}
%\begin{equation}
%\begin{array}{rll}
%\min & b^\top y& \\
%\text{s.t.~} & A^\top y  = \tilde c&\\
%& y_i \ge 0 & i \in [m]\setminus \E \\
%& y_i^+, y_i^- \ge 0 & i \in \E.
%\end{array}\tag{$\tilde D$}\label{tilde D}
%\end{equation}
%%\end{minipage}\par\vspace{\belowdisplayskip}
%\end{prf}
Let \eqref{tilde D}' be obtained from \eqref{tilde D} by replacing each variable $y_i$, $i \in \E$, with $y_i^+ - y_i^-$, where $y_i^+$ and $y_i^-$ are new variables which are required to be nonnegative.
Clearly \eqref{tilde D} and \eqref{tilde D}' are equivalent, so \eqref{tilde D}' has an optimal solution.
Furthermore, since \eqref{tilde D}' is in standard form, it has an optimal solution $\tilde y'$ that is a basic feasible solution.
In particular, via Tardos' algorithm, the vector $\tilde y'$ can be computed in a number of operations polynomially bounded in $\size(A)$.
%\cnote{is this cheating?}
Let $\tilde y$ be obtained from $\tilde y'$ by replacing each pair ${{}\tilde y'_i}^+, {{}\tilde y'_i}^-$ with $\tilde y_i : = {{}\tilde y'_i}^+ - {{}\tilde y'_i}^-$.
It is simple to check that $\tilde y$ is an optimal solution to \eqref{tilde D}.
Since $\tilde y'$ is a basic feasible solution, it has at most $n$ nonzero entries.
By construction, so does $\tilde y$.

%\cnote{Variables $y_i^+, y_i^-$ cannot be both basic in \eqref{tilde D}' since the corresponding columns of the constraint matrix of \eqref{tilde D}' are linearly dependent.}

Next, we discuss how to compute a vector $\tilde y$ that satisfies (i) and (ii).
Let problem \eqref{tilde P}' be obtained from \eqref{tilde P} by dropping the inequalities $a_j^\top x \le b_j$, for $j \in [m] \setminus \E$, such that $a_j$ can be written as a linear combination of $a_i$, $i \in \E$.
Since problem \eqref{tilde P} is feasible, then \eqref{tilde P} and \eqref{tilde P}' have the same feasible region and are therefore equivalent.
Let \eqref{tilde D}' be the dual of \eqref{tilde P}'.
Note that \eqref{tilde D}' is obtained from \eqref{tilde D} by dropping the variables $y_j$ corresponding to the inequalities of \eqref{tilde P} dropped to obtain \eqref{tilde P}'.
Note that \eqref{tilde P}' has the same form of \eqref{tilde P}, thus, from the the first part of the proof, we can compute a vector $\tilde y'$ optimal to \eqref{tilde D}' with at most $n$ nonzero components.
Furthermore, $\tilde y'$ can be computed in a number of operations polynomially bounded in $\size(A)$.
Let $\tilde y$ be obtained from $\tilde y'$ by adding back the dropped components and setting them to zero.
The vector $\tilde y$ is feasible to \eqref{tilde D}, and, from complementary slackness with $\tilde x$, it is optimal to \eqref{tilde D}.
%It is simple to check that $\tilde y$ is optimal to \eqref{tilde D}, and satisfies (i) and (ii).
%From complementary slackness, $\tilde y$ is optimal to \eqref{tilde D}.
Furthermore, $\tilde y$ satisfies (i) and (ii).
% Assume now that $P$ is simple.
% Let \eqref{tilde P}' and \eqref{tilde D}' be defined as in the previous paragraph.
\end{prf}

In the next lemma, we show that at step~\ref{alg: 3} we can always find an index $h \in \H \setminus \E$.

\begin{lemma}
\label{new index}
In step~\ref{alg: 3} of the \iterative, we have $\H \setminus \E \neq \emptyset$.
In particular, the index $h$ exists at each iteration.
%found at  of the  contains at least an index $j \in [m] \setminus \E$.
\end{lemma}

\begin{prf}
Let $\bar c$, $\tilde c$, $\tilde x$ and $\tilde y$ be the vectors computed at a generic iteration of the \iterative. 
Let $\hat c = {\frac{n^3 k \alpha}{\norminf{\bar c}} \bar c}$, and note that $\tilde{c} = \floor{\hat c}$.
Moreover, we have 
$\norminf{\hat c} = n^3 k \alpha$ and, since this number is integer,
% largest absolute value of an entry of $\hat c$ is the integer $n^3 k \alpha$, 
we also have 
$\norminf{\tilde c} = n^3 k \alpha$.

%By construction
Let $\B = \{i \in \{1,\dots,m\} \mid \tilde y_i \neq 0\}$.
From property (i) of the vector $\tilde y$ we know $|\B| \le n$.
%Consider the following partition of $B$:
%\begin{align*}
%B \cap \E \\
%B \setminus \E.
%\end{align*}
From the constraints of \eqref{tilde D} we obtain
\begin{equation}
\label{e: tilde c 1}
\tilde c 
%= B^{\top} \tilde y 
= \sum_{i \in [m]} a_i \tilde y_i
= \sum_{i \in \B} a_i \tilde y_i.
%= \sum_{i \in \B \cap \E} a_i \tilde y_i+ \sum_{i \in \B \setminus \E} a_i \tilde y_i.
\end{equation}
%$B^{\top} \tilde y_B = \tilde c$.

%\smallskip

Note that $\tilde y_j \ge 0$ for every $j \in \B \setminus \E$ since $\tilde y$ is feasible to \eqref{tilde D}.
Hence to prove this lemma we only need to show that 
\begin{align}
\label{killit}
\abs{\tilde y_j} > nk \qquad \text{ for some }j \in \B \setminus \E.
\end{align}
The proof of \eqref{killit} is divided into two cases.

\smallskip

In the first case we assume $\B \cap \E = \emptyset$.
Thus, to prove \eqref{killit}, we only need to show that $\abs{\tilde y_j} > nk$ for some $j \in \B$.
%From the definition of $\B$ we only need to prove $\norminf{\tilde y} \ge nk$.
To obtain a contradiction, we suppose $\abs{\tilde y_j} \le nk$ for every $j \in \B$.
From \eqref{e: tilde c 1} we obtain
\[
\norminf{\tilde c} 
\le \sum_{j \in \B} \norminf{a_j \tilde y_j}
= \sum_{j \in \B} \left(\abs{ \tilde y_j}\norminf{a_j}\right)  
\le \sum_{j \in \B} (nk \cdot \alpha ) 
\le n^2 k \alpha.
%< n^3 k \alpha.
\]
However, this contradicts the fact that $\norminf{\tilde c} = n^3 k \alpha$.
%\note{just need $\norminf{\tilde c} = n^2 k \alpha + 1$ here}. 
Thus $\abs{\tilde y_j} > nk$ for some $j \in \B$, and \eqref{killit} holds.
This concludes the proof in the first case.

\smallskip

%In the second case we assume that $\B \cap \E \neq \emptyset$. In particular, we have $|\B \setminus \E| \le n-1$. To derive a contradiction, suppose that \eqref{killit} does not hold, i.e., $\abs{\tilde y_j} \le nk$ for every $j \in \B \setminus \E$. From \eqref{e: tilde c 1} we obtain 

In the second case we assume that $\B \cap \E$ is nonempty.
In particular, we have $|\B \setminus \E| \le n-1$.
In order to derive a contradiction, suppose that \eqref{killit} does not hold, i.e., $\abs{\tilde y_j} \le nk$ for every $j \in \B \setminus \E$.
%For a contradiction, suppose $\abs{\tilde y_i} < nk$ for $i \in [n]$.
From \eqref{e: tilde c 1} we obtain 
\begin{align*}
%\label{e: tilde c 2}
\tilde c 
%= B^{\top} \tilde y 
= \sum_{i \in \B} a_i \tilde y_i
= \sum_{i \in \B \cap \E} a_i \tilde y_i+ \sum_{j \in \B \setminus \E} a_j \tilde y_j.
\end{align*}
Then
\begin{align}
\label{contra}
\begin{split}
\norminf{\tilde c - \sum_{i \in \B \cap \E} a_i \tilde y_i}
& \le \sum_{j \in \B \setminus \E} \norminf{a_j \tilde y_j}
= \sum_{j \in \B \setminus \E} \left(\abs{ \tilde y_j}\norminf{a_j}\right)  \\
& \le \sum_{j \in \B \setminus \E} (nk \cdot \alpha ) 
\le (n-1) nk \alpha
\le n^2 k \alpha - 1.
\end{split}
\end{align}

Next, in order to derive a 
%In order to complete the proof in this second case, we derive a 
contradiction, we show that 
%But this contradicts Claim~\ref{adpclaim}. % \eqref{adpbound}.
\begin{align}
\label{adpbound}
\norminf{\tilde c - \sum_{i \in \B \cap \E} a_i \tilde y_i}
%\ge \norminf{\hat c -  \sum_{i \in \B \cap \E} a_i \tilde y_i} - \norminf{\tilde c - \hat c}
> n^2 k \alpha - 1.
\end{align}
%\end{claim}
%
%
%$\norminf{\tilde c - \sum_{i \in \B \cap \E} a_i \tilde y_i} \ge ???$.
%
%\note{I would change the chain of implications, and anyway call inequalities when they are needed}
%
%
%
%
%The above claim and \eqref{contra} yield a contradiction, 
%
%
%First, we show the following claim.
%\begin{claim}
%\label{adpclaim}
%
%
%
%\begin{cpf}
By adding and removing $\tilde c$ inside the norm in the left-hand side below, we obtain
% \begin{align*}
% \hat c -  \sum_{i \in \B \cap \E} a_i \tilde y_i
% = \tilde c - \sum_{i \in \B \cap \E} a_i \tilde y_i - (\tilde c - \hat c).
% \end{align*}
% Thus
\begin{align}
\label{thatone}
%\begin{split}
\norminf{\hat c -  \sum_{i \in \B \cap \E} a_i \tilde y_i}
= \norminf{\tilde c - \sum_{i \in \B \cap \E} a_i \tilde y_i - (\tilde c - \hat c)} 
\le \norminf{\tilde c - \sum_{i \in \B \cap \E} a_i \tilde y_i} + \norminf{\tilde c - \hat c}.
%& < \norminf{\tilde c - \sum_{i \in \B \cap \E} a_i \tilde y_i} + 1.
%\end{split}
\end{align}

Let us now focus on the left-hand side of \eqref{thatone}.
We have that $\hat c$ is orthogonal to $a_i$, for every $i \in \E$.
This is because $\hat c$ is a scaling of $\bar c$ and the latter vector is, by definition, orthogonal to $a_i$, for every $i \in \E$.
We obtain
\begin{align}
\label{52}
%\begin{split}
\norminf{\hat c -  \sum_{i \in \B \cap \E} a_i \tilde y_i} 
\ge \frac{1}{ \sqrt n} \norm{\hat c -  \sum_{i \in \B \cap \E} a_i \tilde y_i} 
\ge \frac{\norm{\hat c} }{ \sqrt n} 
\ge \frac{\norminf{\hat c}}{ \sqrt n} 
= \frac{n^3 k \alpha}{ \sqrt n} 
\ge n^2 k \alpha,
%\end{split}
\end{align}
%\note{just need $n^{2.5}k \alpha$ here }
where the second inequality holds by Pythagoras' theorem.

Using \eqref{thatone}, \eqref{52}, and noting that $\norminf{\tilde c - \hat c} < 1$ by definition of $\tilde c$, we obtain
\begin{align*}
\norminf{\tilde c - \sum_{i \in \B \cap \E} a_i \tilde y_i}
\ge \norminf{\hat c -  \sum_{i \in \B \cap \E} a_i \tilde y_i} - \norminf{\tilde c - \hat c}
> n^2 k \alpha - 1.
\end{align*}
This concludes the proof of \eqref{adpbound}.

Inequalities \eqref{contra} and \eqref{adpbound} yield a contradiction, thus \eqref{killit} holds.
This concludes the proof in the second case. %\cnote{fix claim numbering}
%\end{cpf}
%\smallskip
\end{prf}

\subsection{Correctness}

Our next goal is to prove that the \iterative is correct, i.e., that it returns a vertex of $P$ maximizing $c^\top x$.

At each iteration of the \iterative, let $F$ be the face of $P$ defined as 
\[
F := \{ x \in \R^n \mid a_i^\top x \le b_i \text{ for } i \in [m]\setminus \E, \  a_i^\top x = b_i \text{ for } i \in \E\},
\]
and note that $F$ is the feasible region of \eqref{tilde P}.
%
%We will prove that at each iteration the dimension of $F$ decreases by one, and that each optimal solution of \eqref{prob: main} lies in $F$.
%
%The next lemma shows that at each iteration the dimension of $F$ decreases by one. 
The next lemma implies that at each iteration the dimension of $F$ decreases by $1$.
%This will be used to prove that the algorithm performs at most $n$ iterations.
%Note that after $n$ iterations of the algorithm $F$ will be a vertex of $P$.

\begin{lemma}
\label{l: D full row rank}
At each iteration, the row submatrix of $A$ indexed by $\E$ has full row rank.
Furthermore, at each iteration, its number of rows increases by exactly one.
\end{lemma}
%\cnote{change this: the row submatrix of $A$ indexed by $\E$ must have full row rank}

\begin{prf}
We prove this lemma recursively.
Clearly, the statement holds at the beginning of the algorithm because we have $\E = \emptyset$.

Assume now that, at a general iteration, the row submatrix of $A$ indexed by $\E$ has full row rank.
From Lemma~\ref{new index}, the index $h \in \H \setminus \E$ defined in step~\ref{alg: 3} of the algorithm exists.
In particular, $\tilde y_h >nk$.
%the set $\H$ found at step~\ref{alg: 3} contains at least an index $j \notin \E$.
%This means that there exists $h \in [m] \setminus \E$ such that $\tilde y_h \ge  n k$.
From property (ii) of the vector $\tilde y$, we have
%$\tilde y_j = 0$ for every $j \in [m] \setminus \E$ such 
%that $a_j$ can be written as a linear combination of $a_i$, $i \in \E$.
that $a_h$ is linearly independent from the vectors $a_i$, $i \in \E$.
Hence the rank of the row submatrix of $A$ indexed by $\E \cup \{h\}$ is one more than the rank of the row submatrix of $A$ indexed by $\E$.
In particular, it has full row rank.
\end{prf}

In the next three lemmas, we will prove that, at each iteration, every optimal solution to \eqref{prob: main} lies in $F$. Note that, since $F$ is a face of $P$, it is also a $[0,k]$-polytope.

%The complementary slackness conditions state that two vectors $\tilde x$ and $\tilde y$ feasible for \eqref{tilde P} and \eqref{tilde D}, respectively, are optimal if and only if
%\begin{align}
%\tilde y_i >0 \quad \Rightarrow \quad a_i^{\top} \tilde x = b_i && i \in [m]\setminus \E.\label{e: complem}
%\end{align}

Suppose that an optimal solution $\tilde y$ of \eqref{tilde D} is known.
The complementary slackness conditions for linear programming imply that, for every $\tilde x$ optimal for \eqref{tilde P}, we have
\begin{align}
& \tilde y_i >0 \quad \Rightarrow \quad a_i^{\top} \tilde x = b_i && i \in [m]\setminus \E. \label{e: complem}
\end{align}
Thus, in order to solve \eqref{tilde P}, we can restrict the feasible region of \eqref{tilde P} by setting the primal constraints in \eqref{e: complem} to equality.

%\note{Alberto, please check next paragraph.}

Now, suppose that our goal is to solve a variant of \eqref{tilde P} where we replace the cost vector $\tilde c$ with another cost vector $\hat c$. Denote by $(\hat P)$ this new LP problem. Note that $\tilde y$ might not even be feasible for the dual problem $(\hat D)$ associated to $(\hat P)$. Can we, under suitable conditions, still use $\tilde y$ to conclude that a primal constraint is satisfied with equality by each optimal solution to $(\hat P)$?  We show that, if $\tilde y$ is `close' to being feasible for $(\hat D)$, then for each index $i \in [m]\setminus \E$ such that $\tilde y_i$ is sufficiently large, we have that the corresponding primal constraint is active at every optimal solution to $(\hat P)$.
Thus, in order to solve $(\hat P)$, we can restrict the feasible region of $(\hat P)$ by setting these primal constraints to equality.

In the following, for $u \in \R^n$ we denote by $\abs{u}$ the vector whose entries are $\abs{u_i}$, $i=1,\dots n$.
%The next lemma shows (b).

\begin{lemma}
\label{main lemma}
%Let $\textcolor{red}{\tilde P}= \{ x \in \R^n \mid a_i^\top x \le b_i \text{ for } i \in [m]\setminus \E \text{ and } a_i^\top x = b_i \text{ for } i \in \E\}$ be a $[0,k]$-polytope, where $a_i \in \Z^{n}$ and $b_i \in \Z$ for $i \in [m]$.
Let $\tilde x \in F$, $\hat c \in \R^n$, and $\tilde y \in \R^m$ be such that 
\begin{align}
& \abs{A^{\top} \tilde y - \hat c} \le 1 \label{e: quasi feasibility}\\
& \tilde y_i \ge 0 && i \in [m] \setminus \E \label{e: nonnegativity}\\
& \tilde y_i >0 \quad \Rightarrow \quad a_i^{\top} \tilde x = b_i && i \in [m] \setminus \E.\label{e: complementary}
\end{align}
Then for any vector $\hat x \in F \cap \Z^n$ with $\hat c^\top \hat x \ge \hat c^\top \tilde x$, we have %(in particular for any $\hat c$-maximal vector)
\begin{align*}
& \tilde y_i > n k \quad \Rightarrow \quad a_i^{\top} \hat x = b_i && i \in [m] \setminus \E.
\end{align*}
\end{lemma}

\begin{prf}
Let $u := (\hat x - \tilde x)$, and let $u^+, u^- \in \R^n_+$ be defined as follows. For $j \in [n]$,
\begin{align*}
u_j^+ := 
\begin{cases}
u_j & \text{if }u_j \ge 0\\
0& \text{if }u_j< 0,
\end{cases}
\qquad
u_j^- := 
\begin{cases}
0& \text{if }u_j\ge 0\\
-u_j & \text{if }u_j < 0.
\end{cases}
\end{align*}
Clearly $u = u^+ - u^-$ and $\abs{u}=u^+ + u^-$. 
Since $\hat c^\top \hat x \ge \hat c^\top \tilde x$, we have $\hat c ^{\top} u \ge 0$.

We prove this lemma by contradiction.
Suppose that there exists $h \in [m]\setminus \E$ such that $\tilde y_h > n k$ and $a_h^{\top} \hat x \neq b_h$.
Since $\hat x \in F$ and $a_h, \hat x, b_h$ are integral, we have $a_h^{\top} \hat x \le b_h -1$.
We rewrite \eqref{e: quasi feasibility} as
%\begin{align*}
$A^{\top} \tilde y - 1 \le \hat c \le A^{\top} \tilde y + 1$.
%\end{align*}
Thus
\begin{align}
\hat c^{\top} u	&= \hat c^{\top} u^+ - \hat c^{\top} u^- \le (A^{\top} \tilde y + 1)^{\top} u^+ -  (A^{\top} \tilde y - 1)^{\top} u^-\nonumber\\
			& = (A^{\top} \tilde y)^{\top} (u^+ - u^-) +  1^{\top}(u^+ + u^-)= (A^{\top} \tilde y)^{\top} u +  1^{\top}\abs{u}\label{e: 1}.
\end{align}
%\begin{align}
%\hat c^{\top} u	&= \hat c^{\top} u^+ - \hat c^{\top} u^-\nonumber\\
%			&\le (A^{\top} \tilde y + 1)^{\top} u^+ -  (A^{\top} \tilde y - 1)^{\top} u^-\nonumber\\
%			& = (A^{\top} \tilde y)^{\top} (u^+ - u^-) +  1^{\top}(u^+ + u^-)\nonumber\\
%			& = (A^{\top} \tilde y)^{\top} u +  1^{\top}\abs{u}\label{e: 1}.
%\end{align}
We can upper bound $1^{\top}\abs{u}$ in \eqref{e: 1} by observing that $\abs{u_j} \le k$ for all $j\in [n]$, since $u$ is the difference of two vectors in $[0,k]^n$.
Thus
\begin{equation}\label{e: first ub}
1^{\top}\abs{u} \le n k.
\end{equation}
We now compute an upper bound for $(A^{\top} \tilde y)^{\top} u = \tilde y^{\top} A u$ in \eqref{e: 1}.
\begin{align}
\tilde y^{\top} A u
&= \tilde y_h a_h^{\top}u + \sum_{i \in \E}\tilde y_i a_i^{\top}u + \sum_{i \in [m]\setminus \E, \, i \neq h}\tilde y_i a_i^{\top}u \nonumber\\
&<- n k + \sum_{i \in \E}\tilde y_i a_i^{\top}u + \sum_{i \in [m]\setminus \E, \, i \neq h}\tilde y_i a_i^{\top}u \label{e: 4}\\
&=  - n k + \sum_{i \in [m]\setminus \E, \, i \neq h}\tilde y_i a_i^{\top}u \label{e: 2}\\
		&\le - n k + \sum_{i \in [m]\setminus \E, \, i \neq h, \, \tilde y_i>0}\tilde y_i a_i^{\top}u  \label{e: 3}\\
		&\le - n k \label{e: second ub}.
\end{align}
To prove the strict inequality in \eqref{e: 4} we show $\tilde y_h a_h^{\top}u < - n k$.
We have $\tilde y_h > n k >0$, thus condition \eqref{e: complementary} implies $a_h^{\top} \tilde x= b_h$. 
Since $a_h^{\top} \hat x \le  b_h -1$, we get $a_h^{\top} u = a_h^{\top} \hat x - a_h^{\top} \tilde x \le -1$.
We multiply $\tilde y_h > n k$ by $a_h^{\top} u$ and obtain $\tilde y_h \cdot a_h^{\top} u < nk \cdot a_h^{\top} u \le - nk$.
Equality \eqref{e: 2} follows from the fact that, for each $i \in \E$ we have $a_i^{\top}\hat x = b_i$ and $a_i^{\top}\tilde x = b_i$ since both $\hat x$ and $\tilde x$ are in $F$, thus $a_i^{\top}u = 0$.
Inequality \eqref{e: 3} follows from \eqref{e: nonnegativity}.
To see why inequality \eqref{e: second ub} holds, first note that, from condition \eqref{e: complementary}, $\tilde y_i >0$ implies $a_i^{\top}\tilde x = b_i$.
Furthermore, since $\hat x \in F$, we have $a_i^{\top}\hat x \le b_i$.
Hence we have $a_i^{\top} u \le 0$ and so $\tilde y_i a_i^{\top} u \le 0$.

By combining \eqref{e: 1}, \eqref{e: first ub} and \eqref{e: second ub} we obtain $\hat c^{\top} u < 0$.
This is a contradiction since we have previously seen that $\hat c ^{\top} u \ge 0$.
\end{prf}

For a vector $w \in \Z^n$ and a polyhedron $Q \subseteq \R^n$, we say that a vector is \emph{$w$-maximal in $Q$} 
%\note{Do we need this definition? Used 12 times.} 
if it maximizes $w^\top x$ over $Q$.

\begin{lemma}
\label{correctness}
The set $\H$ given at step~\ref{alg: 3} of the \iterative is such that 
every vector $\hat x$ that is $\bar c$-maximal in $F$
%every $\bar c$-maximal vector $\hat x \in F$
satisfies $a_i^{\top} \hat x = b_i$ for every $i \in \H$.
\end{lemma}

\begin{prf}
Clearly, we just need to prove the lemma for every vertex $\hat x$ of $F$ that maximizes $\bar c^\top x$ over $F$. % is $\bar c$-maximal in $F$
In particular, $\hat x$ is a vertex of $P$ and is therefore integral.

Define $\hat c \in \R^n$ as $\hat c_i := {\frac{n^3 k \alpha}{\norminf{\bar c}} \bar c_i}$ for $i = 1,\dots, n$. 
At step~\ref{alg: 3}, 
$\tilde x$ is an optimal vertex of \eqref{tilde P},
%thus it is $\tilde c$-maximal in $F$,
and $\tilde y$ is an optimal solution to the dual \eqref{tilde D}.
%Since $\tilde c = \floor{\hat c}$, 
We have:
\begin{align}
&\abs{A^{\top} \tilde y -  \hat c} \le 1 && \label{e: approx feas}\\
&\tilde y_i \ge 0 && i \in [m] \setminus \E. \label{e: nonneg}
%&\tilde y_i > 0 \Rightarrow a_i^{\top} \tilde x = b_i& i \in [m] \setminus \E.\label{e: compl}
\end{align}
Constraints \eqref{e: nonneg} are satisfied since $\tilde y$ is feasible for \eqref{tilde D}.
Condition \eqref{e: approx feas} holds because $\abs{A^{\top} \tilde y -  \hat c} = \abs{\tilde c -  \hat c} = \hat c - \tilde c < 1$.
%$\hat c -1 \le \floor{\hat c} = \tilde c = \floor{\hat c} \le \hat c +1$. 
Moreover, the complementary slackness conditions 
\eqref{e: complem} 
%\eqref{e: complementary}
are satisfied by $\tilde x$ and $\tilde y$, because they are optimal for \eqref{tilde P} and \eqref{tilde D}, respectively.

Thus, $A,b,\hat c,\tilde x,\tilde y$ satisfy the hypotheses of Lemma \ref{main lemma}.
Since the vector $\hat x$ is $\bar c$-maximal in $F$ and $\hat c$ is a scaling of $\bar c$, the vector $\hat x$ is also $\hat c$-maximal in $F$.
Since $\tilde x \in F$, we have $\hat c^{\top}\hat x \ge \hat c^{\top}\tilde x$.
Then Lemma \ref{main lemma} implies
\begin{align*}
& \tilde y_i > n k  \quad \Rightarrow \quad a_i^{\top} \hat x = b_i && i \in [m]\setminus \E,
\end{align*}
that is, $a_i^{\top} \hat x = b_i$ for all $i \in \H$.
\end{prf}

\begin{lemma}
\label{correctness2}
The set $\E$ updated in step~\ref{alg: 3} of the \iterative is such that
every vector $x^*$ that is $c$-maximal in $P$
%every $c$-maximal vector $x^* \in P$
satisfies $a_i^{\top} x^* = b_i$ for $i \in \E$.
\end{lemma}
\begin{prf}
Consider a vector $x^*$ that is $c$-maximal in $P$.
We prove this lemma recursively.
Clearly, the statement is true at the beginning of the algorithm, when $\E = \emptyset$.

Suppose now that the statement is true at the beginning of a general iteration. % for $\abs{\E} \le h$.
At the beginning of step~\ref{alg: 3} we have that $x^*$ is $c$-maximal in $F$, thus it is also $\bar c$-maximal in $F$. 
When we add an index $h \in \H \setminus \E$ to $\E$ at the end of step~\ref{alg: 3}, by Lemma \ref{correctness} we obtain that $a_h^{\top} x^* = b_h$. Thus, at each iteration of the algorithm we have $a_i^{\top} x^* = b_i$ for $i \in \E$.
\end{prf}

%The \iterative ends if, at step \eqref{alg: 1}, we have $\bar c = 0$. In the next lemma, we show that if this condition is satisfied, then the vector $x^*$ returned by the algorithm solves \eqref{prob: main}.

In the next theorem, we show that the \iterative is correct.

\begin{theorem}
The vector $x^*$ returned by the \iterative 
%at step~\ref{alg: 1} 
is an optimal solution to the LP problem \eqref{prob: main}.
%$\max \{c^\top x \mid x \in P\}$. 
\end{theorem}

%\note{This lemma has no label so it is never used. Strange, it should be!}

\begin{prf}
%Consider a general iteration of the algorithm.
Consider the iteration of the \iterative when the vector $x^*$ is returned at step~\ref{alg: 1}.
Let the set $\E$ be as defined in the algorithm when $x^*$ is returned.
Up to reordering the inequalities defining $P$, we can assume, without loss of generality, that $\E = \{1,\dots,r\}$.
Consider the following primal/dual pair:

\noindent\begin{minipage}{0.5\linewidth}
\begin{equation}
\begin{array}{rll}
\max & c^\top x&\\
\text{s.t.~} & a_i^\top x = b_i & i = 1,\dots,r\\
& a_i^\top x \le b_i & i=r+1,\dots,m
\end{array}\tag{$P$}\label{P}
\end{equation}
\end{minipage}%
\begin{minipage}{0.5\linewidth}
\begin{equation}
\begin{array}{rll}
\min & b^\top y& \\
\text{s.t.~} & A^\top y  = c&\\
& y_i \ge 0 & i=r+1,\dots,m.
\end{array}\tag{$D$}\label{D}
\end{equation}
\end{minipage}\par\vspace{\belowdisplayskip}
Note that the feasible region of \eqref{P} is the same of \eqref{tilde P}, and it consists of the face $F$ of $P$ obtained 
%at each iteration 
by setting to equality all the constraints indexed by $\E$.
Furthermore, the objective function of \eqref{P} coincides with the one of \eqref{prob: main}.

Let $A_{\E}$ be the row submatrix of $A$ indexed by $\E$.
By Lemma \ref{l: D full row rank}, the rank of $A_{\E}$ is $r$.
When, at step~\ref{alg: 1}, we project $c$ onto $\{x \in \R^n \mid  A_{\E}x = 0\}$, we get  $\bar c = c - A_{\E}^\top (A_{\E} A_{\E}^\top)^{-1}A_{\E} c$. 
Since the termination condition is triggered at this iteration, we have $\bar c = 0$.
This implies $c =  A_{\E}^\top z$, where $z := (A_{\E}  A_{\E}^\top)^{-1}A_{\E} c$.

Let $\bar y \in \R^m$ be defined by $\bar y_i := z_i$ for $i=1,\dots,r$, and $\bar y_i := 0$ for $i=r+1,\dots,m$. First, $\bar y$ is feasible for \eqref{D}. In fact $\bar y_i \ge 0$ for $i=r+1,\dots,m$ and
\[
A^\top \bar y = \sum_{i=1}^m \bar y_i a_i
%= \sum_{i=1}^r \bar y_i a_i 
= \sum_{i=1}^r z_i a_i 
= A_{\E}^\top z 
=c.
\]

In particular, $x^*$ is feasible for \eqref{P}.
We have
%Thus, for an arbitrary $x$ feasible for \eqref{P} we have 
\[
c^\top x^* 
= \bar y^\top A x^*
%= z^\top A_{\E} x^*
= \sum_{i=1}^m \bar y_i a_i^{\top} x^*
= \sum_{i=1}^r \bar y_i a_i^{\top} x^* % +  \sum_{i=r+1}^m \bar y_i a_i^{\top} x^* 
= \sum_{i=1}^r \bar y_i b_i
= \bar y^{\top}b.
\] 
By strong duality, $x^*$ is $c$-maximal in $F$. 
If $x^*$ is not $c$-maximal in $P$, then there exist a different vector $x^\dagger$ that is $c$-maximal in $P$. 
In particular, we have $c^{\top}x^\dagger > c^{\top}x^*$.
From Lemma \ref{correctness2}, the vector $x^\dagger$ 
%that is $c$-maximal in $P$ 
lies in $F$.
Since $x^*$ is $c$-maximal in $F$, we obtain $c^{\top}x^\dagger \le c^{\top}x^*$, a contradiction. 
This shows that $x^*$ is $c$-maximal in $P$.
\end{prf}

\subsection{Length of simplex path}

We now present a proof of Theorem~\ref{th 2}, which 
%In next theorem we 
provides a bound on the length of the simplex path generated by the \iterative.

%\begin{theorem}
%\label{th 2}
%The \iterative generates a simplex path from $x^0$ to the output vector.
%The length this simplex path is in $O(n^2 k \log (nk\alpha))$.
%%Given a  vertex $x^0 \in P$ and $c \in \Z^n$, the \iterative generates a simplex path $x^0, x^1,\dots, x^M$ in $P$ such that $x^M$ is an optimal solution to $\max\{c^{\top}x \mid x \in P\}$ and $M \in O(n^6k\log k)$. \cnote{it was $O(n^3k\log k)$}
%\end{theorem}

\begin{prfc}[of Theorem~\ref{th 2}]
First, note that the \iterative constructs a simplex path from $x^0$ to the output vector. 
This is the path obtained by merging all the simplex paths constructed by the \scaling in step~\ref{alg: 3} at each iteration.
From Lemma~\ref{l: feasible}, the \iterative can indeed call the \scaling with input vertex the current $x^*$, since $x^*$ is feasible for problem \eqref{tilde P}.

Next, we upper bound the length of the generated simplex path.
First, note that the \iterative performs at most $n$ iterations.
%we run step~\ref{alg: 3} of the algorithm at most $n$ times. 
This is because, by Lemma \ref{l: D full row rank}, at each iteration the rank of the row submatrix of $A$ indexed by $\E$ increases by one.
Therefore, at iteration $n+1$, the subspace $\{x\in \R^n  \mid a_i^\top x = 0 \text{ for } i \in \E\}$ in step~\ref{alg: 1} is the origin.
Hence the projection $\bar c$ of $c$ onto this subspace is the origin, and the algorithm terminates by returning the current vector $x^*$.

Each time the \iterative performs step~\ref{alg: 3}, it calls the \scaling with input $F$, $x^*$, and $\tilde c$.
Since $F$ is a $[0,k]$-polytope, by Proposition \ref{prop scaling oracle calls}, each time the \scaling is called, it generates a simplex path of length at most $nk (\ceil{\log \norminf{\tilde c}} + 1)$, where $\norminf{\tilde c}= n^3 k \alpha$. % and so $\log \norminf{\tilde c} = 3 \log n + \log k + \log \alpha$.
%
%Denote by $\varphi$ the facet complexity of $P$ and by $\nu$ the vertex complexity of $P$.
%From Theorem~10.2 in \cite{SchBookIP}, we know that facet complexity $\varphi$ and the vertex complexity $\nu$ of $P$ are polynomially related, and in particular $\varphi \le 4 n^2 \nu$. 
%Since $P$ is a $[0,k]$-polytope, we have $\nu \le n \log k$.
%Recall that each inequality in $Ax \le b$ is facet-defining and that the greatest common divisor of the entries in each row of $A$ is one.
%Due to Remark 1.1 in \cite{ConCorZamBook}, we obtain that $\log \alpha \le n \varphi$.
%Hence, we obtain 
%$
%\log \alpha \le n \varphi \le 4 n^3 \nu \le 4 n^4 \log k.
%$ 
%Thus 
%\begin{align*}
%\log \norminf{\tilde c} 
%& = \log(n^3 k \alpha) 
%= 3 \log n + \log k + \log \alpha 
%\le 3 \log n + \log k + 4 n^4 \log k 
%\in O(n^4 \log k).
%\end{align*}
%\cnote{shorten previous formula by removing a step}
Since $\log \norminf{\tilde c} \in O(\log (nk\alpha))$, each time we run the \scaling, we generate a simplex path of length in $O(nk \log (nk\alpha))$.
Therefore, the simplex path generated throughout the entire algorithm has length in $O(n^2 k \log (nk\alpha))$.
\end{prfc}

We immediately obtain the following corollary of Theorem~\ref{th 2}.

\begin{corollary}
If $\alpha$ is polynomially bounded in $n,k$, then the length of the simplex path generated by the \iterative is in $O(n^2 k \log (nk))$.
If we also assume $k=1$, the length reduces to $O(n^2 \log n)$.
\end{corollary}

\subsection{Complexity}

%We now give a bound on the total number of operations performed by the \iterative.
%needed to construct the next vertex in the simplex path.
In this last section, we 
%give a 
bound 
%on 
the number of operations performed to construct the next vertex in the simplex path.

\begin{remark}
\label{rem iterative complexity}
The number of operations performed
%in one iteration of the \iterative
by the \iterative to construct the next vertex in the simplex path is polynomially bounded in $\size(A)$ and $\log k$.
%$n$, $m$, $\log \alpha$, $\log k$.
%If $P$ is simple, this number of operations is in $O(mn^2 + n^5 \log k)$.
If $P$ is simple, the number of operations is in $O(n^4 + n^3 m \log^2 (n k \alpha))$.
\end{remark}

\begin{prf}
%thus the size of $\tilde c$ is at most $n \log \norminf{\tilde c} \in O(n^4 \log k)$.
%Next, we show that $\log \norminf{\tilde c} \in O(n^4 \log k)$ and that $\size(A)$ is at most $n m \log \alpha \in O(m n^5 \log k)$.
%Denote by $\varphi$ the facet complexity of $P$ and by $\nu$ the vertex complexity of $P$.
%From Theorem~10.2 in \cite{SchBookIP}, we know that facet complexity $\varphi$ and the vertex complexity $\nu$ of $P$ are polynomially related, and in particular $\varphi \le 4 n^2 \nu$. 
%Since $P$ is a $[0,k]$-polytope, we have $\nu \le n \log k$.
%Recall that each inequality in $Ax \le b$ is facet-defining and that the greatest common divisor of the entries in each row of $A$ is one \note{We need these assumptions only here? If so, put them in statement here.}.
%Due to Remark 1.1 in \cite{ConCorZamBook}, we obtain that $\log \alpha \le n \varphi$.
%Hence, we obtain 
%$
%\log \alpha \le n \varphi \le 4 n^3 \nu \le 4 n^4 \log k.
%$ 
%Thus 
%\begin{align*}
%\log \norminf{\tilde c} 
%& = \log(n^3 k \alpha) 
%= 3 \log n + \log k + \log \alpha 
%\le 3 \log n + \log k + 4 n^4 \log k 
%\in O(n^4 \log k).
%\end{align*}
%Furthermore, we have  $\log \alpha \le 4 n^4 \log k$, thus 
First, we discuss the number of operations performed in a single iteration of the \iterative:
\begin{enumerate}[label=(\alph*),noitemsep]
\item \label{enum a}
In step~\ref{alg: 1}, computing the projection $\bar c$ of $c$ onto the subspace $\{x\in \R^n  \mid a_i^\top x = 0 \text{ for } i \in \E\}$ can be done in $O(n^3)$ operations via Gaussian elimination.
%in a number of operations polynomial in $n$. % and $\log \alpha$.
%\cnote{we might have $|\E| > n$ even if  $rk A_{\E} \le n$}.
\item 
In step~\ref{alg: 2}, computing the approximation $\tilde c$ of $\bar c$ can be done by binary search, and 
%since $\norminf{\tilde c} = n^3 k \alpha$, 
the number of comparisons required is at most $n\log \norminf{\tilde c}$.
%\item 
%In step~\ref{alg: 3}, at the beginning of each execution of the \basic within the \scaling, we compute an approximation $\tilde c^t$ of $\tilde c$. 
%By construction, $\norminf{\tilde c^t} \le \norminf{\tilde c}$,
%%=n^3 k \alpha$, 
%thus the number of comparisons required is at most $n\log \norminf{\tilde c} \in O(n \log (nk\alpha))$.
%\item 
%The number of operations performed by the \oracle is 
%%in $O(mn)$ if $P$ is simple, due to Remark \ref{rem oracle complexity}.
%%Otherwise, it is 
%polynomially bounded in $\size(A)$.
%and hence
%In the latter case, from the previous discussion on $\size(A)$, it follows that the number of operations performed by the \oracle 
%it is polynomially bounded in $m$, $n$ and $\log k$.
%which is at most $mn \log \alpha$.
%and hence it is polynomially bounded in $m$, $n$ and $\log \alpha$.
\item 
In step~\ref{alg: 3} we call the \scaling to compute the vector $\tilde x$.
From Remark~\ref{rem scaling complexity}, the number of operations performed to construct the next vertex in the simplex path is polynomially bounded in $\size(A)$ and $\log \norminf{\tilde c}$.
If $P$ is simple, the number of operations is $O(n^2m \log^2 \norminf{\tilde c})$.
\item  \label{enum d}
At the end of step~\ref{alg: 3} we compute the vector $\tilde y$.
From Lemma \ref{l: compute dual vars}, the number of operations performed to compute this vector is 
%in $O(n^3)$ if $P$ is simple.
%Otherwise, it is 
polynomially bounded in $\size(A)$, and by $O(nm + n^3)$ if $P$ is simple.
%otherwise. In the latter case, from the previous discussion on $\size(A)$, it follows that the number of operations performed to compute $\tilde y$ is 
%and hence it is polynomially bounded in $m$, $n$ and $\log \alpha$.
\end{enumerate}

Recall from the proof of Theorem \ref{th 2} that the \iterative performs at most $n$ iterations.
Moreover, each vector $\tilde c$ computed at step~\ref{alg: 2} is such that $\log \norminf{\tilde c} \in O(\log (nk\alpha))$.

Let $x^i$ be the current vertex of the simplex path computed by the \iterative, and denote by $x^{i+1}$ the next vertex in the simplex path that we will construct.
%We consider three possible cases.
%\noindent \emph{Case 1}: $x^i$ is the optimal vertex returned by \scaling at step~\ref{alg: 3}.
In the worst case, the vertex $x^{i+1}$ is computed by the \scaling in the very last iteration.
Therefore, the number of operations is bounded by the product of $n$ with the sum of the operations bounds in \ref{enum a}--\ref{enum d} above.
In the general case, this number is polynomially bounded in $\size(A)$ and $\log \norminf{\tilde c}$.
If $P$ is simple, this number is bounded by 
\begin{align*}
O(n \cdot (n^3 + n\log \norminf{\tilde c} + n^2m \log^2 \norminf{\tilde c} + nm + n^3)) 
\in O(n^4 + n^3 m \log^2 \norminf{\tilde c}).
\end{align*}
The statement follows since $\size(A)$ is polynomial in $n,m, \log \alpha$, and 
$\log \norminf{\tilde c} \in O(\log (nk\alpha))$.
\end{prf}

The following remark shows that if the polytope $P$ is `well-described' by the system $Ax \le b$, then the number of operations performed by both the \preprocessingscaling and by the \iterative to construct the next vertex in the simplex path is polynomially bounded in $n$, $m$, $\log k$.
In particular, it is independent on $\alpha$.

\begin{remark}
\label{rem iterative complexity 2}
Assume that $P$ is full-dimensional, that each inequality in $Ax \le b$ is facet-defining, and that the greatest common divisor of the entries in each row of $A$ is one.
Then, the number of operations performed by the \preprocessingscaling and by the \iterative to construct the next vertex in the simplex path is polynomially bounded in $n$, $m$, $\log k$.
\end{remark}

\begin{prf}
From Remark~\ref{rem preprocessing scaling complexity} and Remark~\ref{rem iterative complexity}, the number of operations performed by either algorithm to construct the next vertex in the simplex path is polynomially bounded in $\size(A)$ and $\log k$.
Recall that $\size(A)$ is polynomial in $n$, $m$, and $\log \alpha$.
Therefore, it suffices to show that $\log \alpha$ is polynomially bounded in $n$ and $\log k$.

%Next, we show that $\log \norminf{\tilde c} \in O(n^4 \log k)$ and that $\size(A)$ is at most $n m \log \alpha \in O(m n^5 \log k)$.
Denote by $\varphi$ the facet complexity of $P$ and by $\nu$ the vertex complexity of $P$.
From Theorem~10.2 in \cite{SchBookIP}, we know that $\varphi$ and $\nu$ are polynomially related, and in particular $\varphi \le 4 n^2 \nu$. 
Since $P$ is a $[0,k]$-polytope, we have $\nu \le n \log k$.
%Recall that each inequality in $Ax \le b$ is facet-defining and that the greatest common divisor of the entries in each row of $A$ is one \note{We need these assumptions only here? If so, put them in statement here.}.
Due to the assumptions in the statement of the remark, and Remark 1.1 in \cite{ConCorZamBook}, we obtain that $\log \alpha \le n \varphi$.
Hence, 
$
\log \alpha \le n \varphi \le 4 n^3 \nu \le 4 n^4 \log k.
$ 
%Thus 
%\begin{align*}
%\log \norminf{\tilde c} 
%& = \log(n^3 k \alpha) 
%= 3 \log n + \log k + \log \alpha 
%\le 3 \log n + \log k + 4 n^4 \log k 
%\in O(n^4 \log k).
%\end{align*}
\end{prf}

Note that the assumptions in Remark~\ref{rem iterative complexity 2} are without loss of generality, and it is well-known that we can reduce ourselves to this setting in polynomial time.
Remark~\ref{rem iterative complexity 2} then implies that, if we assume that $k$ is polynomially bounded by $n$ and $m$, then both algorithms run in strongly polynomial time.

To conclude, we highlight that all the obtained bounds on the number of operations performed by our algorithms (see Remarks \ref{rem oracle complexity}--\ref{rem iterative complexity 2}) also depend on the number $m$ of inequalities in the system $Ax \le b$.
This is in contrast with the lengths of the simplex paths, which only depend on $n$ and $k$.
This difference is to be expected, since in order to determine the next vertex, the algorithm needs to read all the inequalities defining the polytope, thus the number of operations must depend also on $m$.

\ifthenelse {\boolean{DCG}}
{
% For DCG begin
\bibliographystyle{spmpsci}
% For DCG end
}
{
% For OO begin
\bibliographystyle{plainurl}
% For OO end
}

\bibliography{../biblio}

\begin{thebibliography}{10}

\bibitem{AckZun95}
D.~Acketa and J.~{\^Z}uni{\'c}.
\newblock On the maximal number of edges of convex digital polygons included
  into an $m \times m$-grid.
\newblock {\em Journal of Combinatorial Theory, Series A}, 69(2):358 -- 368,
  1995.

\bibitem{AhuMagOrl93}
R.K. Ahuja, T.L. Magnanti, and J.B. Orlin.
\newblock {\em Network Flows: Theory, Algorithms, and Applications}.
\newblock Prentice Hall, Englewood Cliffs NJ, 1993.

\bibitem{AloVu97}
N.~Alon and V.H. Vu.
\newblock Anti-hadamard matrices, coin weighing, threshold gates, and
  indecomposable hypergraphs.
\newblock {\em Journal of Combinatorial Theory, Series A}, 79:133--160, 1997.

\bibitem{BalBar91}
A.~Balog and I.~B\'{a}r\'{a}ny.
\newblock On the convex hull of the integer points in a disc.
\newblock In {\em Proceedings of the seventh annual symposium on Computational
  geometry}, SCG '91, pages 162--165, New York, NY, USA, 1991. ACM.

\bibitem{BerTsi97}
D.~Bertsimas and J.~Tsitsiklis.
\newblock {\em Introduction to Linear Optimization}.
\newblock Athena Scientific, Belmont, MA, 1997.

\bibitem{BlaDeLLou20}
M.~Blanchard, J.~De Loera, and Q.~Louveaux.
\newblock On the length of monotone paths in polyhedra.
\newblock {\em arXiv:2001.09575}, 2020.

\bibitem{ConCorZamBook}
M.~Conforti, G.~Cornu\'ejols, and G.~Zambelli.
\newblock {\em Integer Programming}.
\newblock Springer, 2014.

\bibitem{dPMic16}
A.~Del~Pia and C.~Michini.
\newblock On the diameter of lattice polytopes.
\newblock {\em Discrete \& Computational Geometry}, 55(3):681--687, 2016.

\bibitem{DezManOnn18}
A.~Deza, G.~Manoussakis, and S.~Onn.
\newblock Primitive zonotopes.
\newblock {\em Discrete \& Computational Geometry}, 60(1):27--39, July 2018.

\bibitem{DezPou18}
A.~Deza and L.~Pournin.
\newblock Improved bounds on the diameter of lattice polytopes.
\newblock {\em Acta Mathematica Hungarica}, 154(2):457--469, 2018.

\bibitem{DezPouSuk19}
A.~Deza, L.~Pournin, and N.~Sukegawa.
\newblock The diameter of lattice zonotopes.
\newblock {\em Proceedings of the American Mathematical Society (to appear)},
  2020.

\bibitem{Din73}
E.A. Dinic.
\newblock The method of scaling and transportation problems.
\newblock {\em Issled. Diskret. Mat. Science, Moscow. (In Russian.)}, 1973.

\bibitem{EdmKar72}
J.~Edmonds and R.M. Karp.
\newblock Theoretical improvements in algorithmic efficiency for network flow
  problems.
\newblock {\em Journal of the Association for Computing Machinery},
  19:248--264, 1972.

\bibitem{FraTar87}
A.~Frank and \'{E} Tardos.
\newblock An application of simultaneous diophantine approximation in
  combinatorial optimization.
\newblock {\em Combinatorica}, 7:49--65, 1987.

\bibitem{Gab85}
H.N. Gabow.
\newblock Scaling algorithms for network problems.
\newblock {\em Journal of Computer and SystemSciences}, 31:148--168, 1985.

\bibitem{KitMatMiz12}
T.~Kitahara, T.~Matsui, and S.~Mizuno.
\newblock On the number of solutions generated by {D}antzig's simplex method
  for {LP} with bounded variables.
\newblock {\em Pacific Journal of Optimization}, 8(3):447--455, 7 2012.

\bibitem{KitMiz13}
T.~Kitahara and S.~Mizuno.
\newblock A bound for the number of different basic solutions generated by the
  simplex method.
\newblock {\em Mathematical Programming, Series A}, 137:579--586, 2013.

\bibitem{KleOnn92}
P.~Kleinschmidt and S.~Onn.
\newblock On the diameter of convex polytopes.
\newblock {\em Discrete Mathematics}, 102:75--77, 1992.

\bibitem{LenLenLov82}
A.K. Lenstra, H.W. Lenstra, and L.~Lov\'asz.
\newblock Factoring polynomials with rational coefficients.
\newblock {\em Mathematische Annalen}, 261:515--534, 1982.

\bibitem{Miz16}
S.~Mizuno.
\newblock The simplex method using {T}ardos' basic algorithm is strongly
  polynomial for totally unimodular {LP} under nondegeneracy assumption.
\newblock {\em Optimization Methods and Software}, 31(6):1298--1304, 2016.

\bibitem{MizSukDez18}
S.~Mizuno, N.~Sukegawa, and A.~Deza.
\newblock An enhanced primal-simplex based {T}ardos' algorithm for linear
  optimization.
\newblock {\em Journal of the Operations Research Society of Japan},
  61(2):186--196, 2018.

\bibitem{Nad89}
D.J. Naddef.
\newblock The {H}irsch conjecture is true for $(0,1)$-polytopes.
\newblock {\em Mathematical Programming}, 45:109--110, 1989.

\bibitem{SchBookIP}
A.~Schrijver.
\newblock {\em Theory of Linear and Integer Programming}.
\newblock Wiley, Chichester, 1986.

\bibitem{SchBookCO}
A.~Schrijver.
\newblock {\em Combinatorial Optimization. Polyhedra and Efficiency}.
\newblock Springer-Verlag, Berlin, 2003.

\bibitem{SchWeiZie95}
A.S. Schulz, R.~Weismantel, and G.M. Ziegler.
\newblock $0/1$-integer programming: Optimization and augmentation are
  equivalent.
\newblock In {\em Proceedings of ESA '95}, pages 473--483, 1995.

\bibitem{Tar86}
E.~Tardos.
\newblock A strongly polynomial algorithm to solve combinatorial linear
  programs.
\newblock {\em Operations Research}, 34(2):250--256, 1986.

\bibitem{Thi91}
T.~Thiele.
\newblock {E}xtremalprobleme f{\"u}r {P}unktmengen.
\newblock Master's thesis, {F}reie {U}niversit{\"a}t {B}erlin, 1991.

\bibitem{Zie00}
G\"unter~M. Ziegler.
\newblock Lectures on $0/1$-polytopes.
\newblock In Gil Kalai and G\"unter~M. Ziegler, editors, {\em Polytopes ---
  Combinatorics and Computation}, pages 1--41. Birkh\"auser Basel, Basel, 2000.

\end{thebibliography}

\end{document}